\input amstex
\documentstyle{amsppt}
\magnification=\magstep1 \NoRunningHeads

\topmatter

\title
A survey on
spectral multiplicities \\ of ergodic actions
 \endtitle
\author  Alexandre~I.~Danilenko
\endauthor
\abstract
Given a transformation $T$ of a standard measure space $(X,\mu)$, let $\Cal M(T)$ denote the set of spectral multiplicities of the Koopman operator $U_T$ defined in $L^2(X,\mu)\ominus\Bbb C$ by $U_Tf:=f\circ T$.
It is discussed in this survey paper which subsets of $\Bbb N\cup\{\infty\}$ are realizable as $\Cal M(T)$ for various $T$: ergodic, weakly mixing, mixing, Gaussian, Poisson, ergodic infinite measure preserving, etc.
The corresponding constructions are considered in detail.
Generalizations to actions of Abelian locally compact second countable groups are also discussed.

\endabstract

\address
 Institute for Low Temperature Physics
\& Engineering of National Academy of Sciences of Ukraine, 47 Lenin Ave.,
 Kharkov, 61164, UKRAINE
\endaddress
\email alexandre.danilenko\@gmail.com
\endemail

\endtopmatter

\document

\head 1. Introduction
\endhead

In these notes we survey the progress achieved in studying of the {\it spectral multiplicities} of ergodic dynamical systems.
Despite availability of several rather recent nice  sources  on general spectral theory \cite{Go}, \cite{KaTh}, \cite{Le2}, \cite{Le3}, \cite{Na},
 a part of \cite{Ka}, etc.,  the  topic of spectral multiplicities (especially its important developments in the last few years) is not covered there comprehensively.
To fill this gap is the main purpose of the present work.

Let $(X,\goth B,\mu)$ be a standard probability space and let $T$ be an invertible $\mu$-preserving transformation.
Denote by  $U_T$ the induced Koopman unitary operator $U_T$ of $L^2_0(X,\mu):=L^2(X,\mu)\ominus\Bbb C$ given by
$U_Tf:=f\circ T$.
By the spectral theorem for unitary operators
\footnote{
We outline briefly the idea of the proof of the spectral theorem for $U_T$.
Given a real function $f\in L^2_0(X,\mu)$ with $\|f\|_2=1$, the smallest $U_T$-invariant subspace $C(f)\subset L^2_0(X,\mu)$ containing $f$ is called the {\it $U_T$-cyclic subspace generated by $f$}.
Let $\sigma_f$ stand for the probability measure on $T$ such that
$\int_{\Bbb T}z^n\,d\sigma_f(z)=\langle U_T^nf,f\rangle$, $n\in\Bbb Z$.
The linear map sending $U_T^nf$ to the monomial $z^n$ on $\Bbb T$ for each $n\in\Bbb Z$ extends uniquely to a unitary isomorphism of $C(f)$ onto $L^2(\Bbb T,\sigma_f)$ in such a way that $U_T$ corresponds to the operator of multiplication by the independent variable $z\in\Bbb T$.
Next, it is possible to select a sequence $f_1,f_2,\dots$ of real functions in $L^2_0(X,\mu)$ of norm $1$ in such a way that $\bigoplus_j C(f_j)=L^2_0(X,\mu)$ and $\sigma_{f_1}\succeq\sigma_{f_2}\succeq\cdots$.
The equivalence class of $\sigma_{f_1}$ does not depend on the choice of the sequence $(f_i)_i$.
Next, there is a countable Borel partition $S_1,S_2,\dots,S_\infty$ of $\Bbb T$ such that $\sigma_{f_1}\restriction S_j\sim\sigma_{f_j}\restriction S_j$
  and $\sigma_{f_{j+1}}(S_j)=0$ and
$\sigma_{f_1}\restriction S_\infty\sim\sigma_{f_j}\restriction S_\infty$
for each $j<\infty$.
This partition does not depend (by modulo $\sigma_{f_1}$-measure $0$) on the choice of  $(f_i)_i$ either.
Then $L^2_0(X,\mu)$ is unitarily isomorphic to the orthogonal sum
$\bigoplus_j L^2(\Bbb T,\sigma_{f_j})
=\bigoplus_j\bigoplus_{l=1}^jL^2(S_j,\sigma_{f_1})=\int_{\Bbb T}^\oplus\Cal H_z\,d\sigma_{f_1}(z)$, where $\Cal H_z$ is a Hilbert space of dimension $j\in\Bbb N\cup\{\infty\}$ if $z\in S_j$.
The corresponding unitary isomorphism conjugates $U_T$ with the operator of multiplication by the independent variable $z$ on $\Bbb T$.},
 there is a probability measure $\sigma_T$ on $\Bbb T$ and a measurable field of Hilbert spaces $\Bbb T\ni z\mapsto \Cal H_z$ such that
$$
L^2_0(X,\mu)=\int_{z\in\Bbb T}^\oplus\Cal  H_z\,d\sigma_T(z)\quad\text{  and  }\quad U_Tf(z)=zf(z), \ z\in\Bbb T,
$$
for each $f:\Bbb T\ni z\mapsto f(z)\in
\Cal H_z$ with $\int_\Bbb T\|f(z)\|^2\, d\sigma_T(z)<\infty$ \cite{Nai}.
The measure $\sigma_T$ is called a {\it measure of maximal spectral type} of $U_T$.
Consider a map $m_T:\Bbb Z\ni z\mapsto m_T(z)=\text{dim}\,\Cal H_z\in\Bbb N\cup\{\infty\}$.
It is called the {\it spectral multiplicity map} of $T$.
Denote by $\Cal M(T)$ the essential range of $m_T$ (with respect to $\sigma_T$)\footnote{In a similar way  the set $\Cal M(V)$ is defined for an arbitrary unitary operator  $V$ in a separable Hilbert space.}.
We can now state the {\it spectral multiplicity problem} as follows:
\roster
\item"{\bf (Pr1)}"
Given a subset $E\subset\Bbb N\cup\{\infty\}$, is there an ergodic transformation $T$ such that $\Cal M(T)=E$?
\endroster
If the answer is affirmative, we say that $E$ is {\it realizable}.
We consider (Pr1) as a weak version of one of  the main problems in the spectral theory of dynamical systems:
\roster
\item"{\bf (Pr2)}" Which unitary operator with continuous spectrum is unitarily equivalent to a Koopman operator?
\endroster
The two problems are  open.
Currently, neither restrictions on the realizable sets nor on Koopman realizable operators (except  for the obvious ones) are known.
However in contrast with the extremely difficult~(Pr2), there is a significant progress in tackling (Pr1).
The list of known realizable sets is impressive.
It includes
\roster
\item"---" all subsets $E\ni 1$,
\item"---" all subsets $E\ni 2$,
\item"---" the  subsets $n\cdot E:=\{ne\mid e\in E\}$ for arbitrary sets $E\ni 1$ and $n>2$,
\item"---"  all multiplicative (and additive) subsemigroups  of $\Bbb N$,  \item"---" the subsets $\{\#(\Gamma\backslash\goth S_n/\goth S_k)\mid k=1,\dots,n-1\}$, where $\Gamma$ is a subgroup in the full symmetric group $\goth S_n$, the inclusion $\goth S_k\subset \goth S_n$ is standard, and the symbol $\#$ denotes the cardinality,
    \item"---" the sets
    $\{n,m,nm\}$, $\{l,m,n, lm,ln,mn, lmn\}$, \dots.
\endroster
To show this,
  several subtle constructions were elaborated.
 We present here  a modern detailed  exposition of  them.
We first note that it suffices to consider (Pr1) only in the case when $\infty\not\in E$.
The point is that we can always join $\infty$ via a simple trick.
Indeed, if  $T$ has a singular spectrum (this is always the case for the constructions considered below) then $\Cal M(T\times B)=\Cal M(T)\cup\{\infty\}$ for a Bernoullian transformation $B$.
The following  three main technologies are used to attack (Pr1): isometric extensions, Cartesian products and auxiliary group actions.
Sections 2, 3 and 5 of these notes are devoted to them respectively.

\subhead (Tec1) Isometric extensions
\endsubhead
The main idea of this method was invented by Oseledets,
who constructed the first example of an ergodic transformation with non-simple spectrum of finite multiplicity \cite{Os}.
For that he considered  a {\it double Abelian} or {\it meta-Abelian} extension of an interval exchange maps.
 The smallest  possible non-commutative meta-Abelian group
$(\Bbb Z/2\Bbb Z)\ltimes(\Bbb Z/3\Bbb Z)$ was the {\it fiber} in his extensions.
The structure of semidirect product gives rise to a specific kind of {\it symmetry} of the extension.
This symmetry generates non-trivial values of the spectral multiplicity map of the extension.
The paper \cite{Os} had a strong influence on subsequent progress in solving (Pr1).
Robinson generalized Oseledets example by considering more general metacyclic
$(\Bbb Z/p\Bbb Z)\ltimes(\Bbb Z/q\Bbb Z)$-extensions of rank-one transformations.
He  showed  that for each integer $n>0$, the set $\{1,n\}$ is realizable \cite{Ro}.
In a subsequent paper \cite{Ro} he proved that, more generally,
each finite subset $E\subset\Bbb N$ satisfying the two conditions
\roster
\item"(C1)" $1\in E$ and
\item"(C2)" if $n,m\in E$ then lcm$(n,m)\in E$
\endroster
is realizable.
Robinson's argument consists of two parts.
First, he  solves  certain group theoretical problem (we call it an {\it algebraic realization}): to put in correspondence to $E$ some metacyclic group, say $K$.
Then he applies the powerful Katok-Stepin periodic approximation technique \cite{KaSt} to  prove the existence of  an appropriate  $K$-extension via generic (Baire category) argument.
We note that only finite subsets can be realized this way.
In order to realize infinite $E$ satisfying (C1) and (C2), Goodson, Kwiatkowsky, Lema\'nczyk and Liardet modified Oseledets-Robinson construction.
They replaced $K$ with infinite compact Abelian groups \cite{Go-Li}.
Since this modification destroyed the symmetry inherent to the double extensions, they imposed  an additional  symmetry condition on their {\it single} Abelian extensions.
In a subsequent paper \cite{KwiLe} the condition (C2) was removed by passing to natural factors of the Abelian extensions satisfying  this symmetry condition.
Thus, it was proved that each subset $E\subset\Bbb N$  containing $1$ is realizable (even in the class of weakly mixing transformations).
In \cite{Da4} this result was reproved in the framework
of double extensions, i.e. in the spirit of \cite{Os}.
The fibers of the extensions considered in \cite{Da4} are homogeneous spaces of some special compact meta-Abelian groups.
We note that the realizations from \cite{Go-Li}, \cite{KwiLe} and \cite{Da4} are constructed in an explicit effective (cutting-and-stacking) way.
Some versions of the algebraic realization problem play an important role in the three papers.
Other ergodic and weakly mixing realizations of $E\subset\Bbb N$, $1\in E$, appeared in \cite{Ag3} and \cite{Ag5} respectively.
At the first glance, they look different from the Oseledets' double extensions and do not involve any algebraic realizations.
We show, however, that Ageev's constructions are based indeed on the same ideas.

\subhead  (Tec2) Cartesian products
\endsubhead
Katok suggested to consider Cartesian products as a source of ergodic transformations with non-trivial spectral multiplicities (see \cite{Ka} which circulated since the mid-eighties as an unpublished manuscript).
The main interest was to realize subsets not containing $1$.
He showed, in particular, that for a generic transformation $T$, $\Cal M(T\times T)\subset\{2,4\}$ and conjectured that $\Cal M(T^{\times n})=\{n,n(n-1),\dots,n!\}$ for each $n>1$.
This conjecture was proved simultaneously by Ryzhikov  \cite{Ry2} (the particular case $n=2$) and Ageev \cite{Ag2} (the general case).
Concrete examples of such transformations are discussed in \cite{Ag4},
\cite{Ry3}, \cite{Ry4}, \cite{PrRy}.
They include the classical Chacon map, some Adams staircase transformations \cite{Ad}, del Junco-Rudolph rank-one map \cite{dJRu}.

\subhead (Tec3) Auxiliary group actions
\endsubhead
We consider this approach to spectral multiplicities as a natural development of the method of Cartesian products because it utilizes the same kind of symmetries.
 However it  is free of the restrictions on the spectral measure which are inherent to Cartesian products.
 Ageev first used  auxiliary group actions in \cite{Ag6} to answer affirmatively the following Rokhlin's question\footnote{In the case $n=2$ the answer was given earlier in \cite{Ag2} and \cite{Ry2}.} which is a particular case of (Pr1):
\roster
\item"{\bf (Pr3)}" given $n>1$, is  $\{n\}$ is realizable?
 \endroster
The main idea  is to select an appropriate meta-Abelian countable group $G_n$ and fix an element $g\in G_n$ in such a way that
for a generic action $T=(T_g)_{g\in G_n}$ of $G_n$, the transformation $T_g$ is ergodic and $\Cal M(T_g)=\{n\}$.
Ageev's argument was further developed and simplified in \cite{Da1} and \cite{Ry5}.
An explicit construction of weakly mixing transformations with homogeneous spectrum of an arbitrary multiplicity appeared in \cite{Da1}.

Subsequent progress was connected   with various {\bf combinations} of these 3 techniques.
By  considering  isometric extensions of transformations with homogeneous spectrum, the present author constructed weakly mixing realizations of  $n\cdot E$, where $n>1$ and $E$ is an arbitrary subset of $\Bbb N$ containing $1$ \cite{Da1}.
For that we combined~(Tec1) and~(Tec3).
This is explained in Section~5.
On the other hand, Katok and Lema{\'n}czyk  considered isometric extensions of Cartesian squares of weakly mixing transformations and constructed realizations of all finite subsets  of $\Bbb N$ containing $2$ \cite{KaLe}.
\footnote{In fact, an extra condition on the realizable subsets was imposed in \cite{KaLe}. However, as was shown later in \cite{Da4}, this condition is satisfied for all finite subsets of $\Bbb N$.}
Extending their result, we show in \cite{Da3} that every  subset of $\Bbb N$ containing $2$ is realizable.
A combination of  (Tec1) and (Tec2) is used in \cite{KaLe} and \cite{Da3}. This is explained in detail in Section~4.
As was shown in \cite{Ry5}, a combination of (Tec2) and (Tec3) also leads to interesting results.
For instance, for all $n,m\in\Bbb N$, the subset $\{n,m,nm\}$ is realizable.
We discuss this in Section~5.

Further in these notes we study (Pr1) under some restrictions.
In Section~6,  (Pr1) is considered in the class of {\it mixing} transformations.
The first examples of mixing transformations with  simple spectrum were found in the class of  Gaussian transformations \cite{New}.
More elementary, rank-one\footnote{Gaussian transformations are never of  rank one \cite{dR}.},  systems appeared in
  \cite{Or} (stochastic constructions) and \cite{Ad} (concrete examples).
For each $n>1$, Robinson constructed in \cite{Ro2} a mixing $T$ with
$n<\max\Cal M(T)<\infty$.
He also asked about existence of a mixing $T$ with $\max\Cal M(T)=n$.
Ageev answered this affirmatively in \cite{Ag7}.
We conjecture that if  $E=\Cal M(T)$ for an ergodic transformation $T$ then there is also a mixing $T$ with the same property.
At present, this is verified for all $E$ realizable via (Tec1), (Tec2)
and their combination.
The existence of mixing solution of (Pr3) (with usage of (Tec3)) is proved in a recent work \cite{Ti2}.
The main idea is to obtain mixing realizations as a limit of non-mixing ones.
It originates from \cite{Ry3} (see also \cite{Ry4}).
As was shown in \cite{Ag7} and \cite{Ry4}, there is a mixing transformation $T$ such that all the symmetric powers $T^{\odot n}$, $n\in\Bbb N$, have simple spectrum.
This led to mixing realizations of  the following sets  $\{\#(\Gamma\backslash\goth S_n/\goth S_k)\mid k=1,\dots,n-1\}$ of cardinalities of double cosets, where $\Gamma$ is a subgroup of the symmetric group $\goth S_n$.
In \cite{Da4} we construct mixing realizations for each $E\subset\Bbb N$ that contains  $1$ or $2$.

In Section 7 we consider (Pr1) in the class of infinite measure preserving transformations.
Since $L^2(X,\mu)$ does not contain constants
 if $\mu(X)=\infty$, the Koopman operator $U_T$
is considered in the entire $L^2(X,\mu)$.
Quite surprisingly,  (Pr1) can be solved {\it completely} in this case by combining (Tec1) and (Tec2) \cite{DaRy1}.
Moreover, this result was refined in a subsequent paper \cite{DaRy2}: for each $E\subset\Bbb N$, there exists a mixing (or zero type) ergodic multiply recurrent infinite measure preserving transformation\footnote{In the infinite case, mixing does not imply ergodicity and ergodicity does not imply recurrence.} $T$ with $\Cal M(T)=E$.
Thus, we conclude that the main obstacle to realizations of spectral multiplicities  in the probability preserving case is the presence of constants in $L^2(X,\mu)$.

Section 8 is devoted to Gaussian and Poisson realizations.
Girsanov  showed  in \cite{Gi} that for an ergodic Gaussian transformation $T$, either $T$ has a simple spectrum or $\Cal M(T)$  is infinite.
By \cite{DaRy2}, every subsemigroup of $\Bbb N$ (either additive or multiplicative) admit a mixing Poisson (and hence Gaussian) realization.
Some other (non-semigroup) subsets, for example $\{1, 1\cdot 3,1\cdot 3\cdot 5,\dots\}$ are also Poisson realizable.

In Section 9 we briefly survey the state of the spectral multiplicity problem for actions of general locally compact second countable Abelian groups \cite{DaLe}, \cite{DaSo}, \cite{LePa}.
Section 10 contains a list of open problems and some additional material on the spectral multiplicities.
Appendix is of algebraic nature.
It is devoted completely to the algebraic realizations which play a role in Sections~2, 4 and 9.

The author thanks the anonymous referee for the useful remarks which improved and corrected the original version of these notes.

\head
2. Isometric extensions.
Realization of subsets containing $1$
\endhead

\subhead Isometric extensions
\endsubhead
Let $T$ be a transformation with a simple spectrum on the space $(X,\goth B,\mu)$.
Let $K$ be a compact second countable group.
Denote by $\lambda_K$ the Haar measure on it.
Let $\phi:X\to K$ be a Borel map.
We will call it a {\it cocycle} of $T$.
Then we can consider a new transformation $T_\phi$ on the product space
$(X\times K,\mu\times\lambda_K)$:
$$
T_\phi(x,k):=(Tx,k\phi(x)).
$$
It is called  a {\it compact group  (or skew product) extension} of $T$.
We note that $(T_\phi)^n=(T^n)_{\phi^{(n)}}$, where
$$
\phi^{(n)}(x):=\cases
 \phi(x)\phi(Tx)\cdots\phi(T^{n-1}x) &\text{if }n>0,\\
1&\text{if }n=0.
\endcases
$$
Two cocycles $\phi,\psi:X\to K$ are called {\it cohomologous} (we  write $\phi\asymp\psi$) if there is a function $a:X\to K$ such that $\phi(x)=a(x)\psi(x)a(Tx)^{-1}$.
Then the corresponding extensions $T_\phi$ and $T_\psi$ are isomorphic.
Given  a closed subgroup $Y$ of $K$, we consider a transformation
$T_{Y\backslash\phi}$ of the space $(X\times Y\backslash K,\mu\times \lambda_{Y\backslash K})$:
$$
T_{Y\backslash\phi}(x,Yk):=(Tx,Yk\phi(x)),
$$
where $\lambda_{Y\backslash K}$ is the Haar measure on the homogeneous space $Y\backslash K$.
This transformation is  called an {\it isometric  extension} of $T$.

\subhead
Main result and a strategy of its  proof
\endsubhead
We now state the main result of this section.

\proclaim{Theorem 2.1} For each subset $E\subset\Bbb N$, there is a weakly mixing transformation $R$
such that $\Cal M(R)=E\cup\{1\}$.
\endproclaim

The desired transformation $R$  is constructed as an isometric extension of another transformation with a simple spectrum.
We present here two different  (but close) constructions.

Suppose that $K$ is Abelian.
Denote by $\widehat K$ the dual group to $K$.
Since $K$ is compact and second countable, $\widehat K$ is countable and discrete.
There is a natural decomposition of $L^2_0(X\times K,\mu\times\lambda_K)$ into an orthogonal sum
of $U_{T_\phi}$-invariant subspaces
 $$
L^2_0(X\times K,\mu\times\lambda_K):=L_0^2(X,\mu)\oplus\bigoplus_{1\ne\chi\in\widehat K}
L^2(X,\mu)\otimes\chi.
$$
The restriction of $U_{T_\phi}$ to the subspace $L^2(X,\mu)\otimes\chi$
is unitarily equivalent to the unitary operator $U_{T_\phi,\chi}$ acting in $L^2(X,\mu)$ by
$$
(U_{T_\phi,\chi}f)(x):= \chi(\phi( x)) f(Tx).
$$
Thus we obtain a decomposition of $U_{T_\phi}$ into an orthogonal sum
$$
U_{T_\phi}=\bigoplus_{\chi\in\widehat K}U_{T_\phi,\chi},\tag2-1
$$
where $U_{T_\phi,1}$ denotes $U_T$.
Since there is a canonical embedding of the dual group $\widehat {Y\backslash K}$ into $\widehat K$, we will always consider
$\widehat {Y\backslash K}$ as a subgroup of $\widehat K$.
Then  \thetag{2-1} yields also a decomposition for $U_{T_{Y\backslash\phi}}$:
$$
U_{T_{Y\backslash\phi}}=\bigoplus_{\chi\in\widehat{Y\backslash K}}U_{T_\phi,\chi}.\tag2-2
$$
Denote by $\sigma_{T,\phi,\chi}$ a measure of maximal spectral  type of $U_{T_\phi,\chi}$.
Suppose now that the following properties hold.
\roster
\item"(P1)"
$U_{T_\phi,\chi}$ has a simple spectrum for each $\chi$.
\item"(P2)" There is an equivalence relation $\approx$ on $\widehat K$
such that
\item"(P3)"
$\sigma_{T,\phi,\chi}\sim\sigma_{T,\phi,\chi'}$ if $\chi\approx\chi'$ and
\item"(P4)"
$\sigma_{T,\phi,\chi}\perp\sigma_{T,\phi,\chi'}$ if $\chi\not\approx\chi'$.
\endroster
Then it follows from \thetag{2-2} that
$$
\Cal M(T_{Y\backslash\phi})=\{\#(\overset{\centerdot}\to\chi\cap \widehat{Y\backslash K})\mid 0\ne \chi\in \widehat{Y\backslash K}\},
\tag2-3
 $$
where $\overset{\centerdot}\to\chi$ denotes the $\approx$-equivalence class of $\chi$.

\subhead Simple facts about unitary operators
\endsubhead
Given a unitary operator $U$, we denote by WCP$(U)$ the weak closure
of the powers of $U$.
Every $U$-invariant subspace is invariant under each operator from the semigroup WCP$(U)$.
If \rom{WCP}$(U)\ni\alpha I$ with $|\alpha|<1$ then the maximal spectral type of  $U$ is continuous.
The next lemma follows directly from the spectral theorem for unitary operators (see also \cite{Go-Li, Proposition~3}).
It is our main tool   to establish (P4).

\proclaim{Lemma 2.2} Let $U$ and $V$ be two unitary operators.
If \rom{WCP}$(U\oplus V)\ni \alpha I\oplus
\beta I$ and $\alpha\ne \beta$ then the  maximal spectral types of $U$ and $V$ are orthogonal.\footnote{This lemma implies, in particular, that if $T$ is $\kappa$-weakly mixing, i.e. WCP$(U_T)\ni(1-\kappa)U_T$, and $0<\kappa<1$ then the convolutions of the maximal spectral type of $T$ are pairwise disjoint \cite{St}.}
\endproclaim

The following well known lemma (see, e.g. \cite{Na, Theorem~4.6}) will be repeatedly used below.

\proclaim{Lemma 2.3} Let $U$ be a unitary operator in the Hilbert space $\Cal H$.
Let $(\Cal H_n)_{n=1}^\infty$ be a sequence of $U$-cyclic subspaces in $\Cal H$.
Suppose that $\lim_{n\to\infty}\inf_{v\in\Cal H_n}\|h-v\|=0$ for each $h\in\Cal H$.
Then $U$ has a simple spectrum.
\endproclaim

\subhead Rank-one transformations
\endsubhead
To satisfy (P1) we consider rank-one transformations.
A transformation $T$ of $(X,\goth B,\mu)$ is called of {\it rank one} if it admits a sequence of {\it Rokhlin towers} generating the entire $\sigma$-algebra
$\goth B$.
This means that there is a sequence of subsets $(A_n)_{n=1}^\infty$ and a sequence of positive integers $(h_n)_{n=1}^\infty$
such that $T^iA_n\cap T^jA_n=\emptyset$ whenever  $0\le i\ne j<h_n$
and
$$
\lim_{n\to\infty}\min_{J\subset\{0,1,\dots,h_n-1\}}\mu(A\triangle \bigsqcup_{j\in J}T^jA_n)=0
$$
for each Borel subset $A\subset X$.
Every rank-one transformation is ergodic.
It follows from Lemma~2.3 that every rank-one transformation
   has a simple spectrum \cite{Ba}.
   The converse is not true \cite{dJ2}.
Every ergodic rotation on a compact group is of rank one \cite{dJ1}.
Every rank-one transformation can be constructed via a simple geometric {\it cutting-and-stacking} inductive process as follows.
Suppose we are given a sequence $(r_n)_{n=1}^\infty$ of positive integers $r_n\ge 2$ and a sequence $(s_n)_{n=1}^\infty$ of maps $s_n:\{0,1,\dots,r_n-1\}\to\Bbb Z_+$.
We construct inductively a sequence of {\it towers}, i.e. finite collections of intervals, called {\it levels}, of the same length.
 The levels are thought as being placed one above the other.
Suppose that the $n$-th tower is of {\it hight $h_n$}, i.e. it consists of $h_n$ levels $I_n^0,\dots, I_n^{h_n-1}$
numbered  from bottom to top.
We cut this tower into $r_n$ subtowers of equal width and enumerate them from left to right.
Then each level $I_n^i$ is partitioned into $r_n$ sublevels
which we denote by $I_n^i(0), \dots,I_n^i(r_n-1)$.
Put now $s_n(i)$ additional levels, called {\it spacers}, over
the $i$-th subtower.
Stack these {\it extended} subtowers by placing the $i$-th subtower on the top of $(i-1)$-th one.
We thus obtain a new tower of hight $h_{n+1}=r_nh_n+\sum_{i=0}^{r_n-1}s_n(i)$.
It is $r_n$ times thinner then the $n$-th tower.
Renumbering the levels of this tower  from bottom to top as $I_{n+1}^0,\dots, I_{n+1}^{h_{n+1}-1}$
we obtain the $(n+1)$-th tower of the construction.
Now we define a transformation $T$ on the union of the first $h_n-1$ levels of the $n$-th tower
 as a translation of each such level one level above.
When $n\to\infty$, we obtain a measure preserving transformation $T$ of a $\sigma$-finite measure space $(X,\goth B,\mu)$.
The measure $\mu$ is finite if and only if $\sum_{n=1}^\infty\frac{\sum_{i=0}^{r_n-1}s_n(i)}{h_n}<\infty$.
It is easy to verify that $T$ is  of rank-one.
We call it the {\it rank-one transformation associated with $(r_n, s_n)_{n=1}^\infty$}.
We note that there are plenty of sequences $(r_n,s_n)_{n=1}^\infty$ which determine the very same (up to isomorphism) rank-one transformation.

By choosing $(r_n,s_n)_{n=1}^\infty$ in an appropriate way we can manufacture various weak limits of powers of $T$.

\proclaim{Lemma 2.4 \rom{(see, e.g., \cite{Ry4})}}
Let  $r_{n_k}\to\infty$
for some  $n_k\to \infty$ .
\roster
\item"\rom{(i)}"
If  $s_{n_k}\equiv 0$ then $U_T^{ph_{n_k}}\to I$ for each integer $p>0$.
\item"\rom{(ii)}"
If  $s_{n_k}\equiv 1$ then $U_T^{ph_{n_k}}\to U_T^{-p}$ for each integer $p>0$.
\item"\rom{(iii)}"
If  $s_{n_k}(i)=i$ for all $i=0,\dots,r_{n_k}-1$
then $U_T^{ph_{n_k}}\to 0$ for each integer $p>0$.
\item"\rom{(iv)}"
If
$$
s_{n_k}=
\cases
0 &\text{\rom{if }}0\le i<r_{n_k}/2,\\
1 &\text{\rom{if }}r_{n_k}/2\le i<r_{n_k}
\endcases
$$
\endroster
then $U_T^{ph_{n_k}}\to 0.5(I+U_T^{-p})$ for each integer $p>0$.
\endproclaim

 In a similar way we can construct  rank-one systems with other {\it polynomial} weak limits of powers.
In this connection it is interesting to note the following  {\it universal} property of the weak limits  from Lemma~2.4(iv).

\proclaim{Proposition 2.5 \cite{Ry4}}
If $0.5(I+U_T^{-k})\in \text{\rom{WCP}}(U_T)$ for an ergodic transformation $T$ and each integer $k>0$
then $p(U_T)\in\text{\rom{WCP}}(U_T)$
 for each polynomial $p$ with non-negative coefficients and such that $p(1)\le 1$.
\endproclaim

\subhead Cocycles of product type
\endsubhead
We are going to isolate an important class of cocycles of rank-one systems.
Such cocycles will be used in the proof of Theorem~2.1.
Let $T$ be a rank-one transformation associated with a sequence
$(r_n,s_n)_{n=1}^\infty$.
Suppose we are given a sequence $(a_n)_{n=1}^\infty$ of maps $a_n:\{0,1,\dots,r_n-1\}\to K$ with $a_n(0)=1$.
We now define inductively a sequence of maps $\beta_n:X_n\to K$, where
$X_n$ is the union of all levels of the $n$-th tower.
We set
$\beta_1(x)=1$ for all $x\in X_1$ and
$$
\beta_{n+1}(x):=\cases
\beta_n(x)a_n(j) &\text{if }x\in I_n^i(j)\text{ for some $0\le i<h_n$ and $0\le j<r_n$},\\
1&\text{if $x\in X_{n+1}\setminus X_n$}
\endcases
$$
for $n>0$.
We note that $\beta_{n}$ is constant on every level of the $n$-th tower.
Hence we can think of it as a {\it labeling} of the levels of the $n$-th tower with elements of $K$.
We now define a cocycle $\phi:X\to K$ of $T$ by setting
$$
\phi(x)=\beta_n(x)\beta_n(Tx)^{-1}\qquad\text{if }x\in X_n\setminus I_n^{h_n-1}.\tag2-4
$$
It is easy to verify that $\phi$ is well defined.
We note that $\phi$  is constant on every level (except of the highest one) of each tower.
We call $\phi$ a {\it cocycle of product type of $T$ associated with $(a_n)_{n=1}^\infty$}.
It is also called a {\it Morse cocycle} (see \cite{Go} and references therein).
Cocycles of product type are convenient to model various weak limits of the corresponding compact extensions of the underlying rank-one maps.

Suppose that $r_{n_l}\to\infty$,  $s_{n_l}\equiv 0$.
It follows from \thetag{2-4} and the definition of $\beta_{n_l+1}$ that
$$
\phi^{(h_{n_l})}(x)=
\beta_{n_l+1}(x)\beta_{n_l+1}(T^{h_{n_l}}x)^{-1}
=\beta_{n_l}(x)a_{n_l}(j)a_{n_l}(j+1)^{-1}\beta_{n_l}(x)^{-1}\tag2-5
$$
for $x\in I^{i}_{n_l}(j)$, $0\le j<r_l-1$, $0\le i<h_{n_l}-1$.
We introduce some notation.
Given $k\in K$, we let $\widetilde k(x,h):=(x,kh)$ for $(x,h)\in X\times K$.
Then $\widetilde k$ is a transformation of $(X\times K,\mu\times\lambda_K)$ commuting with $T_\phi$.

\proclaim{Lemma 2.6} If $K$ is Abelian\footnote{To emphasize the fact that $K$ is Abelian we will use the  ``additive'' symbol $+$ for the group operation in $K$.},
$r_{n_l}\to\infty$,  $s_{n_l}\equiv 0$ and $a_{n_l}(j)-a_{n_l}(j+1)=k$ for all $0\le j<r_{n_l}-1$ then $(U_{T_\phi})^{h_{n_l}}\to U_{\widetilde k}$ as $l\to\infty$.
More generally, given $k_1,\dots,k_p\in K$, if
$$
\frac{\#\{j\mid a_{n_l}(j)-a_{n_l}(j+1)=k_s\}}{r_{n_l}}\to\frac 1p
$$
for each $s=1,\dots,p$ then $(U_{T_\phi})^{h_{n_l}}\to p^{-1}\sum_{s=1}^pU_{\widetilde{k_s}}$.
\endproclaim
\demo{Proof}
To show the first claim of the lemma it is enough to note that
${(T_\phi)}^{h_{n_l}}=(T^{h_{n_l}})_{\phi^{(h_{n_l})}}$,
$U_{T}^{h_{n_l}}\to I$ weakly by Lemma~2.4(i) and
$\phi^{(h_{n_l})}(\mu)\to
\delta_k$ in the $*$-weak topology in view of \thetag{2-5}.
Here $\delta_k$ stands for the Dirac measure supported at $k$.
The second claim is proved in a similar way.
\qed
\enddemo

\subhead Proof of Theorem~2.1 via meta-Abelian isometric extensions
\endsubhead
We outline here a slightly simplified version of the proof given in \cite{Da4}.
Without loss of generality we may assume that
$\emptyset\ne E\not \ni 1$.
 Let $D$ be a compact monothetic totally disconnected group,  $G$ a discrete countable (topological) $D$-module and $H$ a subgroup in $G$.
Given $\chi\in H$, we denote by $\Cal O(\chi)$ the $D$-orbit of $\chi$.
Let
$$
L(G,D,H):=\{\#(\Cal O(\chi)\cap H)\mid 0\ne \chi\in H\}.
$$
By Lemma~A.5, there is a triplet $(G,D,H)$ such that $L(G,D,H)=E$.
Let $K$ denote the dual group to $G$.
Since the duality  preserves the module structure, $K$ is a topological $D$-module.
Let $Z:=D\ltimes K$.
Then $Z$ is a meta-Abelian compact group.
Recall that the multiplication in $Z$ is given by
$$
(d,k)(d',k'):=(d+d',k+d\cdot k'), \qquad d,d'\in D,k,k'\in K.
$$
We identify canonically the dual group $\widehat{G/H}$ to $G/H$ with a closed subgroup, say $Y$, in $K$.
We note that  $\widehat{K/Y}=H$.

Our purpose  now is to construct a rank-one transformation $T$ on a probability space $(X,\goth B,\mu)$ and a product-type-cocycle $\tau:X\to Z$ of $T$ such that
the isometric extension $T_{Y\backslash\tau}$ of $T$ is weakly mixing and $\Cal M(T_{Y\backslash\tau})=E\cup\{1\}$.
For that we will write $T_{Y\backslash\tau}$ as an Abelian extension of a system with a simple spectrum to
satisfy the conditions (P1)--(P4) and make use of \thetag{2-3}.
Thus from now on let $T$ be  associated
with a sequence $(r_n,s_n)_{n=1}^\infty$ and let
 $\tau$ be associated with a sequence $(a_n)_{n=1}^\infty$.
Since $Z$ is the product $D\times K$ as a topological space,  we consider  the cocycle $\tau$  as a pair $(\psi,\phi)$ of maps $\psi:X\to D$ and $\phi:X\to K$.
We now have
$$
T_{Y\backslash\tau}(x,d,k+Y)=(Tx,d+\psi(x),k+d\cdot\phi(x)+Y).
$$
Hence we can think of $T_{Y\backslash\tau}$  as a {\it double} compact Abelian extension, i.e. a compact Abelian extension $(T_\psi)_{Y\backslash\widetilde\phi}$ of a compact Abelian extension
$T_\psi$ of $T$,
where $\widetilde\phi:X\times D\ni (x,d)\mapsto d\cdot\phi(x)\in K$ is a cocycle of $T_\psi$.
Next, for each $n$, the map $a_n:\{0,\dots,r_n-1\}\to Z$ is a pair $(d_n,b_n)$ of maps $d_n:\{0,\dots,r_n-1\}\to D$
and $b_n:\{0,\dots,r_n-1\}\to K$.
We note that $\psi$ is a cocycle of product type associated with
the sequence $(d_n)_{n=1}^\infty$.
In the following claim we explain how to satisfy (P1).

\proclaim{Claim 2.7}
Let $d$ be an element of $D$ generating a dense subgroup in $D$.
Let $n_l\to\infty$.
If $r_{n_l}\to\infty$, $s_{n_l}\equiv 0$,
$b_{n_l}\equiv 0$ and $d_{n_l}(j)-d_{n_l}(j+1)=d$ for all $0\le j<r_{n_l}-1$
then $T_\psi$ is of rank one and
for each $ \chi\in G$, the unitary operator
$U_{T_\tau,\chi}:=U_{(T_\psi)_{\widetilde\phi},\chi}$ has a simple spectrum.
\endproclaim
\demo{Proof}
Since $D$ is compact and $G$ is countable, the $D$-orbit of $\chi$ is finite.
Therefore there is
a nested sequence  $D_1\supset D_2\supset\cdots$ of open subgroups in $D$
such that $\bigcap_{j=1}^\infty D_j=\{0\}$ and $d'\cdot\chi=\chi$ for all $d'\in D_1$.
Then we obtain an increasing  sequence of  subspaces
$$
L^2(X\times D/D_1)\subset L^2(X\times D/D_2)\subset\cdots
$$
in $L^2(X\times D)$ whose union is dense in $L^2(X\times D)$.
Every such a subspace is invariant under $U_{T_\tau,\chi}$.
Hence in view of Lemma~2.3 it is enough to prove that the restriction of $U_{T_\tau,\chi}$ to
$L^2(X\times D/D_j)$ has a simple spectrum for each $j$.
Notice that since $D$ is monothetic, the quotient $D/D_j$ is a finite
cyclic group.
Let
$$
K_\chi:=\{k\in K\mid \chi(d'\cdot k)=1\text{ for all }d'\in D\}.
$$
Then $K/K_\chi$ is a finite\footnote{Indeed, $K_\chi=\bigcap_{d'\in D}\text{Ker}(d'\cdot\chi)$, the $D$-orbit of $\chi$ is finite and the subgroup $\text{Ker}(d'\cdot\chi)$ is of finite index in $K$.
Hence $K_\chi$ is also of finite index in $K$.}
 $D/D_j$-module for each $j$.
Thus we see that it suffices to prove the claim only in the particular case when
$K$ is finite and $D$ is cyclic.
It follows from \thetag{2-5} and the conditions of the claim that
$$
\tau^{(h_{n_l})}(x)=\beta_{n_l}(x)(d,0)\beta_{n_l}(x)^{-1}= (d, c_l(x)-d\cdot c_l(x))
$$
for all  $x\in I_{n_l}^i(j)$, $0\le i<h_{n_l}$ and $0\le j<r_{n_l}-1$,
where $c_l(x)\in K$ is the second coordinate of $\beta_{{n_l}}(x)$.
Since
\roster
\item"---"
$
U_{T_\tau,\chi}^{h_{n_l}}F(x,d')=\chi\circ d'(c_l(x)-d\cdot c_l(x))F(T^{h_{n_l}}x,d+d')
$
for $F\in L^2(X\times D)$, $d'\in D$ and $x$ as above and
\item"---"
 the function $X_{n_l}\ni x\mapsto c_l(x)\in K$ is constant on each level of the $n_l$-th
tower,
\endroster
the $U_{T_\tau,\chi}$-cyclic subspace generated by the vector $F:=1_{I_{n_l}^0}\otimes 1_{\{0\}}$ {\it almost} contains vectors
$1_{I_{n_l}^i}\otimes 1_{\{d'\}}$ for all $0\le i<h_{n_l}-1$ and $d'\in D$.
Hence by Lemma~2.3, $U_{T_\tau,\chi}$ has a simple spectrum.

Substituting $\chi=1$ in the above reasoning, we obtain that $T_\psi$ is of rank one.
\qed
\enddemo

Now we introduce an equivalence relation $\approx$ on $\widehat{K}=G$ to satisfy~(P2).
We call two characters  $\approx$-{\it equivalent} if they belong to the same $D$-orbit.
Then~(P3) is easy to verify.
Indeed,  it is straightforward that
$$
\widetilde\phi\circ\widetilde d=
d\cdot\widetilde\phi\text{ \ for each }d\in D.\tag2-6
$$
This implies that the unitary operator $U_{\widetilde d}$
intertwines $U_{T_\tau,\chi}$ with $U_{T_\tau,d\cdot\chi}$.
This yields (P3).

We now explain how to satisfy (P4).
We first note that by Lemma~A.5, the $D$-module $K$ is finitary, i.e.
 there is a dense countable subgroup $\Cal K$ in $K$
such that the $D$-orbit $\Cal O_*(k)$ of $k$ is finite for each $k\in\Cal K$.
Given $\chi\in H$, we let $l_\chi:=\frac 1{\# \Cal O(\chi)}\sum_{\eta\in\Cal O(\chi)}\eta\in L^2(K)$.

\proclaim{Claim 2.8}
Let $k\in\Cal K$ and $n_l\to\infty$.
If $r_{n_l}\to\infty$, $s_{n_l}\equiv 0$,
$d_{n_l}\equiv 0$
and for each $q\in \Cal O_*(k)$,
$$
\frac{\#\{j\mid b_{n_l}(j)-b_{n_l}(j+1)=q\}}{r_{n_l}}\to\frac 1{\#\Cal O_*(k)}
\tag2-7
$$
then
$U_{T_\tau,\chi}^{h_{n_l}}\to l_\chi(k)I$ for each $\chi\in\widehat K$  as $l\to\infty$.
\endproclaim

\demo{Proof}
We  first note that
$$
\frac 1{\#\Cal O_*(k)}\sum_{q\in \Cal O_*(k)}\eta(q)=\int_D\eta(d\cdot k)\,d\lambda_D(d)=\int_D(d\cdot\eta)(k)\,d\lambda_D(d)=l_\eta(k)
\tag2-8
$$
 for each $\eta\in G$, where $\lambda_D$ stand for the normed Haar measure on $D$.

It follows from the conditions of the claim  and \thetag{2-5} that
$$
\tau^{(h_{n_l})}(x)=\beta_{n_l}(x)(0,b_{n_l}(j)-b_{n_l}(j+1))
\beta_{n_l}(x)^{-1}= (0, c_l'(x)\cdot (b_{n_l}(j)-b_{n_l}(j+1)))
$$
for all  $x\in I_{n_l}^i(j)$, $0\le i<h_{n_l}$ and $0\le j<r_{n_l}-1$,
where $c_l'(x)\in D$ is the first coordinate of $\beta_{{n_l}}(x)$.
Therefore
$$
U_{T_\tau,\chi}^{h_{n_l}}F(x,d')
=\chi((d'+c_l'(x))\cdot (b_{n_l}(j+1)-b_{n_l}(j)))F(T^{h_{n_l}}x,d')
\tag2-9
$$
for each $F\in L^2(X\times D)$, $d'\in D$ and $x$ as above.
It follows from \thetag{2-7}  that the mapping
$$
I^i_{n_l}\supset I^i_{n_l}(j)\ni x\mapsto b_{n_l}(j+1)-b_{n_l}(j)\in K
$$
maps the measure $\mu\restriction I^i_{n_l}$ to a measure which is close (uniformly in $i$) to the uniform distribution, say $\vartheta$, on
$\Cal O_*(k)\subset K$ when $l$ is large.
We now note that the map $X\ni x\mapsto c_l'(x)\in D$
is constant on the levels $I^i_{n_l}$ of the $n_l$-th
tower.
Hence the image of $\mu$ under the mapping
$$
\bigsqcup_i I^i_{n_l}\supset I^i_{n_l}(j)\ni x\mapsto (d'+c_l'(x))(b_{n_l}(j+1)-b_{n_l}(j))\in K
$$
tends to  $\vartheta$ uniformly in $d'\in D$ as $l\to\infty$.
Hence passing to the limit in \thetag{2-9} and using the facts that
$\int_K\chi\,d\vartheta=l_\chi(k)$ by \thetag{2-8}
and $T^{h_{n_l}}\to\text{Id}$ by Lemma~2.4(i)
 we obtain that
$U_{T_\tau,\chi}^{n_l}\to l_\chi(k)I$  as $l\to\infty$.
\qed
\enddemo

Now if $\chi\not\approx\chi'$ then $l_\chi\perp l_{\chi'}$ in $L^2(K)$.
Hence there is $k\in \Cal K$ such that $l_\chi(k)\ne l_{\chi'}(k)$.
Suppose that Claim~2.8 holds for this $k$.
Then Lemma~2.2 yields that the measures of maximal spectral type of
$U_{T_\tau,\chi}$ and $U_{T_\tau,\chi'}$ are orthogonal, i.e. (P4) is satisfied.

In the next claim we show how to make $T_{Y\backslash\tau}$ weakly mixing.

\proclaim{Claim 2.9}
Suppose that $T_\psi$ is ergodic and $l_\chi(k)I\in\text{\rom{WCP}}(U_{T_\tau,\chi})$ for all $\chi\in G$
and $k\in\Cal K$.
If $r_{n_l}\to\infty$, $s_{n_l}(j)=0$ if $0\le j<r_{n_l}/2$ and
$s_{n_l}(j)=1$ if $r_{n_l}/2\le j<r_{n_l}$ and  $d_{n_l}\equiv 0$
for some subsequence $n_l\to\infty$
then $T_{Y\backslash\tau}$ is weakly mixing.
\endproclaim
\demo{Proof}
Since $E\not \ni 1$, it follows that $\#\Cal O(\chi)>1$ for each $\chi\ne 1$.
Therefore there is
 $k\in\Cal K$ with $|l_\chi(k)|<1$.
 Since $l_\chi(k)I\in\text{WCP}(U_{T_\tau,\chi})$, the maximal spectral type of $U_{T_\tau,\chi}$ is continuous.
If $\chi=1$ then $U_{T_\tau,\chi}=U_{T_\psi}$.
 It follows from the conditions of the claim and Lemma~2.4(iv)  that $U_T^{h_{n_l}}\to 0.5(I+U_T^*)$.
Since $d_{n_l}\equiv 0$, we obtain that
$U_{T_\psi}^{h_{n_l}}\to 0.5(I+U_{T_\psi}^*)$.
Therefore the only possible eigenvalue for
$U_{T_\psi}$ is $1$.
However $1$ can not be an eigenvalue of $U_{T_\psi}$ because $T_\psi$ is ergodic (it is of rank one by Claim~2.7).
Thus $U_{T_\psi}$ has a continuous spectrum as well.
Since $U_{T_{Y\backslash\tau}}=\bigoplus_{\chi\in H}U_{T_\tau,\chi}$,
it follows that $T_{Y\backslash\tau}$ is weakly mixing.
\qed
\enddemo

To summarize, we  start with a  sequence $r_n\to\infty$.
Then we partition $\Bbb N$ into infinitely many  subsequences $n_l\to\infty$.
On each of these subsequences  we define $s_{n_l}$ and $a_{n_l}$ in an appropriate way described in  Claims 2.7--2.9.
Then we obtain a rank-one transformation $T$ associated with $(r_n,s_n)_{n=1}^\infty$ and a product-type cocycle $\tau$ of $T$
with values in $K$ such that the transformation $T_{Y\backslash\tau}$ is weakly mixing and the corresponding conditions (P1)--(P4) are satisfied
for the extension $(T_\psi)_{Y\backslash\widetilde\phi}=T_{Y\backslash\tau}$.
Then $\Cal M(T_{Y\backslash\tau})=L(G,D,H)\cup\{1\}=E\cup\{1\}$.
The proof of Theorem~2.1 is complete.

\subhead Generic properties of ``double'' Abelian extensions
\endsubhead
The group Aut$(X,\mu)$ of all measure preserving transformations of a
standard probability space $(X,\mu)$ is Polish in the {\it weak topology} which is induced from the weak operator topology on the group  $\Cal U(L^2(X,\mu))$ of unitary operators in $L^2(X,\mu)$ via the embedding
$T\mapsto U_T$.
For each compact second countable group $Z$, the  set $\Cal A(X,Z)$
of measurable maps from $X$ to $Z$ is Polish when endowed with
the topology of convergence in measure.
If $T$ is aperiodic then it follows from the
Rokhlin lemma that
\roster
\item"(i)"
the conjugacy class $\{RTR^{-1}\mid R\in \text{Aut}(X,\mu)\}$ of $T$
is dense in Aut$(X,\mu)$ (see, e.g., \cite{Na}),
\item"(ii)"
for each cocycle $\phi\in\Cal A(X,Z)$ of $T$, the cohomology class
$\{\psi\in\Cal A(X,Z)\mid\psi\asymp\phi\}$ of $\phi$
is dense in $\Cal A(X,Z)$ \cite{Sc}.
\endroster
If $Z$ is Abelian and $\chi$ is a continuous character of $Z$ then the map
$$
\text{Aut}(X,\mu)\times\Cal A(X,Z)\ni (T,\phi)\mapsto U_{T_\phi,\chi}\in\Cal U(L^2(X,\mu))
$$
is continuous.
We recall a couple of well known facts from the spectral theory of unitary operators.
For the proof we refer to \cite{Na}.
There is a continuous map  $\sigma:\Cal U(L^2(X,\mu))\ni U\mapsto\sigma(U)$
from  $\Cal U(L^2(X,\mu))$ to the Polish space of probability measures on $\Bbb T$ endowed with the $*$-weak topology such that $\sigma(U)$ is a measure of maximal spectral type of $U$.
The set of unitary operators with a simple continuous spectrum is a dense $G_\delta$ in $\Cal U(L^2(X,\mu))$.

Now let $E, G, Z, Y$ be as in the proof of Theorem~2.1.

We let $\Gamma:=\text{Aut}(X,\mu)\times \Cal A(X,Z)$.

\proclaim{Theorem 2.10}
The subset
$$
\Cal S:=\{(T,\tau)\in\Gamma \mid
T_{Y\backslash\tau}\text{ is weakly mixing and } \Cal M (T_{Y\backslash\tau})=E\cup\{1\}\}
$$
is residual in $\Gamma$.
\endproclaim
\demo{Proof}
For all $\chi,\chi'\in G$, the subsets
$$
\align
M_\chi &:=\{(T,\tau)\in\Gamma\mid U_{T_\tau,\chi}
\text{ has a simple continuous spectrum}\},\\
L_{\chi,\chi'}&:=\{(T,\tau)\in\Gamma\mid \sigma(U_{T_\tau,\chi})\perp \sigma(U_{T_\tau,\chi'})\}
\endalign
$$
are $G_\delta$ in $\Gamma$.
Hence the intersection
$$
\Cal J:=\bigcap_{\chi\in G}M_\chi
\cap\bigcap_{\chi\not\approx\chi'\in G}L_{\chi,\chi'}
$$
is also a $G_\delta$ in $\Gamma$.
Of course, $\Cal J\subset \Cal S$.
Therefore it suffices to show that
$\Cal J$ is dense in $\Gamma$.
We introduce  an equivalence relation on $\Gamma$ by
setting
$(T,\tau)\sim(T',\tau')$ if $R^{-1}TR=T'$ and $\tau\circ R\asymp\tau'$.
It is routine to verify that if $(T,\tau)\in\Cal J$ and
$(T',\tau')\sim(T,\tau)$ then $(T',\tau')\in\Cal J$.
It follows from Theorem~2.1 that $\Cal J$ is non-empty.
The properties (i) and (ii) imply that the equivalence class of each pair
$(T,\tau)\in \Cal J$ is dense in $\Gamma$.
\qed
\enddemo

\remark{Remark \rom{2.11}}
In earlier works \cite{Ro1} and \cite{Ro2} Robinson used a category argument to prove the existence of realizations for some subsets with $1$ (see also \cite{Ag3} and \cite{Ag5} for a development of this idea).
We think however that  in  the modern ergodic theory it is easier  to construct first some concrete realizations (as in the previous subsections) and then use them in a standard generic argument to show
that the desired realizations are residual in certain topological space of parameters.
We also note that if $Z$ is a compact Abelian group then the set
$\{(T,\phi)\mid \Cal M(T_\phi)=\{1\}\}$
is residual in Aut$(X,\mu)\times\Cal A(X,Z)$ \cite{Ro4}.
\endremark

\subhead Remark on Ageev's approach from \cite{Ag3} and \cite{Ag5}
\endsubhead
In \cite{Ag3} and \cite{Ag5} Ageev gives  alternative  proofs of Theorem~2.1.
 The realizations from \cite{Ag3} have a discrete component in the spectrum
 and  the realizations from \cite{Ag5} are weakly mixing.
Ageev considers some special single compact Abelian extensions over   rank-one maps.
We will show that these extensions are, in fact, double Abelian (as in the proof of Theorem~2.1).
An additional extension is ``hidden'' in the base  which  is itself a compact Abelian extension of another rank-one map.
Let $E=\{1,n_1,n_2,\dots\}$ and let $F=\bigotimes_k(\Bbb Z/2\Bbb Z)^{n_k}$.
Ageev introduces the following  equivalence relation $\cong$ on the dual group
$\widehat F=\bigoplus_{k}(\Bbb Z/2\Bbb Z)^{n_k}$:
two characters $\chi=(\chi_{1,1},\dots,\chi_{1,n_1},\chi_{2,1}, \dots,\chi_{2,n_2},\dots)$ and
$\chi'=(\chi_{1,1}',\dots,\chi_{1,n_1}',\chi_{2,1}', \dots,\chi_{2,n_2}',\dots)$ are $\cong$-equivalent if
either $\chi=\chi'$ or there exist $k$ and integers  $i,j$ such that $1\le i\ne j\le n_{k}$
such that $\chi_{k,i}$ and $\chi_{k,j}'$ are the only non-trivial coordinates of $\chi$ and $\chi'$ respectively.
We now explain how it is related to algebraic realizations.
It was shown in \cite{KwiLe} that for each $k$, there exist $m_k>n_k$ and a subset $A_k\subset\{1,\dots,m_k\}$ with
$\# A_k=n_k$ such that if
$$
\align
G  &:=\bigoplus_{k}(\Bbb Z/2\Bbb Z)^{m_k},\\
H &:=\{\chi=(\chi_{1,1},\dots,\chi_{1,n_1},\chi_{2,1}, \dots,\chi_{2,n_2},\dots)\in G\mid \chi_{i,k}=1\text{ if }i\not\in A_k\,\forall k\},\\
v &:=\bigoplus_{k} v_k,
\endalign
$$
where $v_k$ is the automorphism of $(\Bbb Z/2\Bbb Z)^{m_k}$ generated by
the cyclic permutation of coordinates,
then
$\cong$ is the $v$-orbit equivalence relation restricted
to $H$ (we identify $\widehat F$ with $H$ in a natural way).
Hence $E=L(G,v,H)$ (see Appendix for the definition of this set).
Passing to a compactification  as in the proof of Lemma~A.5 we obtain a compact monothetic group $D$ acting on $G$  with the same orbits as $v$.
Therefore $E=L(G,D,H)$.
We recall a standard notation.
Given a transformation $R\in\text{Aut}(X,\mu)$, the group of $\mu$-preserving transformations commuting with $R$ is called the {\it centralizer} of $R$.
It is denoted by $C(R)$.

It follows from the proof of Theorem~2.1 that
there is  a rank-one probability preserving transformation $T$ and  a cocycle $\tau=(\psi,\phi)$ of $T$ with values in $D\ltimes K$
such that $T_{Y\backslash\tau}$ is weakly mixing and $\Cal M(T_{Y\backslash\tau})=E$, where $K:=\widehat G$ and $Y:=\widehat{G/H}$.
Moreover,  $T_\psi$ is of rank one  according to Claim~2.7.
The cocycle $Y\backslash\widetilde\phi$ of $T_\psi$ takes values in
the group
$Y\backslash K=F$.
Since  $C(T_\psi)\supset\{\widetilde d\mid d\in D\}$ and $\widetilde\phi\circ \widetilde d=d\cdot\widetilde\phi$,
it follows that if $\chi\cong\chi'\in \widehat F$
then there is a transformation $S\in C(T_\psi)$ such that
$\chi'\circ\widetilde\phi=\chi\circ\widetilde\phi\circ S$.
To summarize,  there exist a rank-one transformation
$Q:=T_\psi$ of a probability space $(Z,\kappa)$,
 measurable maps $\beta_1,\beta_2,\dots$ from $Z$ to $\Bbb Z/2\Bbb Z$ and transformations $S_{i,j}\in C(Q)$, $1\le j< n_i$ such that
if we set
$$
\beta(z):=(\beta_1(z),\beta_1(S_{1,1}z),\dots, \beta_1( S_{1,n_1-1}z),
\beta_2(z),\beta_2(S_{2,1}z),\dots)\in F,
$$
$z\in Z$, then $\beta=Y\backslash\widetilde\phi$ and hence the compact group extension $Q_\beta$ of $Q$ is weakly mixing and $\Cal M(Q_\beta)=E$.
 Such realization of $E$ was obtained in \cite{Ag5} via category argument.

\subhead Proof of Theorem 2.1 via Abelian  extensions (\cite{Go-Li}, \cite{KwiLe})
\endsubhead
By Lem\-ma~A.3,
there exist a countable group $G$, a subgroup $H$ of $G$ and an automorphism $v$ of $G$ such that $E\cup\{1\}=L(G,v,H)$.
We let $K:=\widehat G$, $Y:=\widehat{G/H}$.
We consider a rank-one transformation $T$ associated  with a sequence $(r_n,s_n)_{n=1}^\infty$ and a product-type cocycle $\phi$ of $T$ with values in $K$.
Let $\phi$ be associated with a sequence $(b_n)_{n=1}^\infty$.
Our purpose is to choose the parameters $(r_n,s_n,b_n)_{n=1}^\infty$
in such a way that the skew product transformation $T_\phi$  is weakly mixing and the conditions  (P1)--(P4) are all satisfied.
The condition (P1) holds automatically.
We denote by $\approx$ the $\widehat v$-orbit equivalence relation on $K$, where $\widehat v$ is the automorphism of $K$ dual to $v$.
To satisfy (P3) we need a ``symmetry'' condition on $\phi$ instead of  \thetag{2-6} which is no longer available.
For that we modify  Claim~2.8 as follows.

\proclaim{Claim 2.12}
Let $\Cal K$ be a dense countable subgroup in $K$ such that the $\widehat v$-orbit of every point of $\Cal K$ is finite.
If for each $k\in\Cal K$, there are  two sequences $n_l\to\infty$ and $z_l\to\infty$
such that $r_{n_l}\to\infty$, $s_{n_l}\equiv 0$, \thetag{2-7} holds,
 $z_{n_l}/r_{n_l}\to 0$ very fast and $b_{n_l}(j+z_l)=v(b_{n_l}(j))$ for all
$0\le j<r_n-z_n$.
Then the sequence $T^{z_1h_{n_1}+\cdots+z_lh_{n_l}}$ converges to a transformation $S\in C(T)$ as $l\to\infty$.
 Moreover, $U_{T_\phi}^{h_{n_l}}\to l_\chi(k)I$ and $\phi\circ S\approx v\circ\phi$.
\endproclaim
The latter condition implies (P3) and the former condition yields (P4).
Hence by \thetag{2-3}, $\Cal M(T_{Y\backslash\phi})=E\cup\{1\}$, as desired.

\head 3. Cartesian products
\endhead

\subhead Realization of the set $\{2\}$
\endsubhead
We start this section with two auxiliary lemmata about unitary operators.

\proclaim{Lemma 3.1} Let $V$ be a unitary operator with a simple continuous spectrum in a Hilbert space $\Cal H$.
If $\text{\rom{WCP}}(V)\ni aI+bV$ for some  $a,b\in\Bbb C\setminus\{0\}$ then  $\Cal M(V\otimes V)=\{2\}$.
Moreover, the unitary $W$ of $\Cal H\otimes\Cal H$ given by $W(h_1\otimes h_2):=Vh_2\otimes h_1$,
has a simple spectrum and $W^2=V\otimes V$.
\endproclaim

\demo{Proof}
The proof of the first claim consists of two steps.

(A) Let $h$ be a cyclic vector for $V$.
Denote by $\Cal C$ the smallest $(V\otimes V)$-invariant subspace
containing $h\otimes h$ and $h\otimes Vh$.
Since $aI+bV\in\text{\rom{WCP}}(V)$, it follows that $\Cal C$ is invariant under $(aI+bV)\otimes (aI+bV)$.
Hence $\Cal C$ is also invariant under  the operator
$
I\otimes V+ V\otimes I.
$
This implies that the vectors $h\otimes V^nh$, $n\in\Bbb Z$, are all in $\Cal C$.
Hence  the vectors $V^mh\otimes V^nh$ for arbitrary  $n,m\in\Bbb Z$ also belong to $\Cal C$.
Therefore $\Cal C=\Cal H\otimes\Cal H$
and hence $\Cal M(V\otimes V)\le 2$.

(B)
Let $\sigma$ denote a measure of maximal spectral type of $V$.
Then $\sigma\times\sigma$ is a measure of maximal spectral type
of the unitary representation
$(n,m)\mapsto V^n\otimes V^m$ of $\Bbb Z^2$ in $\Cal H\otimes\Cal H$.
Let
$$
\sigma\times\sigma=\int_\Bbb T\sigma_z\,d\sigma^{*2}(z)
$$
stand for the
 disintegration of $\sigma\times\sigma$
 with respect to the projection map $\Bbb T\times\Bbb T\ni (z_1,z_2)\mapsto z_1z_2\in\Bbb T$.
Then  the convolution  square $\sigma^{*2}$ of $\sigma$ is a measure of maximal spectral type of $V\otimes V$
and the map $\Bbb T\ni z\mapsto \text{dim}(L^2(\Bbb T^2,\sigma_z))$
is the multiplicity function of $V\otimes V$.
Since $\sigma$ is continuous,
the  $\sigma\times\sigma$-measure of the diagonal is zero.
Since this measure is invariant under the flip $(z_1,z_2)\mapsto (z_2,z_1)$, it follows that
$\sigma_z$ is invariant under the flip for a.a. $z$.
Hence if $\sigma_z$ is not continuous then the number of atoms of $\sigma_z$ is even.
It follows that  the values of the multiplicity map of $V\otimes V$ are either even or infinity.

Now (A) plus (B) yield $\Cal M(V\otimes V)=\{2\}$.
As for the second claim, it is enough to note that $W(h\otimes Vh)=Vh\otimes Vh$ and $W^2=V\otimes V$.
Hence $h\otimes h$ belongs to the $W$-cyclic subspace generated by $h\otimes Vh$.
Hence $h\otimes Vh$ is a $W$-cyclic vector in $\Cal H\otimes\Cal H$.
\qed
\enddemo

The first claim of the  following lemma is a generalization of Lemma~2.2.
It follows directly from the spectral theorem for unitary operators.
The second claim was established in \cite{Ry2}.

\proclaim{Lemma 3.2}
Let $V$ and $W$ be  unitary operators with continuous spectrum.
Let $\sigma_V$ and $\sigma_W$
denote measures of maximal spectral type of $V$ and $W$.
\roster
\item"\rom{(i)}"
If there are polynomials or, more generally, analytic functions  $p\ne q$ such that \rom{WCP}$(V\oplus W)\ni p(V)\oplus q(W)$ then $\sigma_V\perp\sigma_W$.
\item"\rom{(ii)}"
 If
\rom{WCP}$(V)$ contains $\alpha I+(1-\alpha)V$ for some  $0<\alpha<1$ then $\sigma_V\perp(\sigma_V)^{*2}$.
\endroster
\endproclaim

The following  is a partial solution of (Pr3):
 the set $\{2\}$ is realizable.

\proclaim{Theorem 3.3  (\cite{Ag2}, \cite{Ry2}, \cite{Ry3})}
Let $T$ be a weakly mixing transformation
on $(X,\goth B,\mu)$ with a simple spectrum.
If $\alpha I+(1-\alpha)U_T\in\text{\rom{WCP}}(U_T)$ for some $0<\alpha<1$ then  $\Cal M(T\times T)=\{2\}$.
Moreover, the transformation $S$ of $(X\times X,\mu\times\mu)$, given by $S(x,y)=(Ty,x)$,
has a simple spectrum and $S^2=T\times T$.
\endproclaim
\demo{Proof}
We note that
$U_{T\times T}$ is unitarily equivalent to the orthogonal sum of $U_T\otimes U_T$ and two copies of $U_T$.
It remains to  apply Lemmata~3.1 and 3.2(ii).
\qed
\enddemo

\subhead Spectral multiplicities of higher Cartesian products
\endsubhead
Given a unitary operator $U$ in the Hilbert space $\Cal H$ and a subgroup $\Gamma$ in the symmetric group $\goth S_n$, we denote by $U^{\otimes n/\Gamma}$ the restriction of the the unitary operator $U^{\otimes n}$
to the subspace of $\Gamma$-invariant tensors in $\Cal H^{\otimes n}$.
In particular,  $U^{\otimes n/\goth S_{n}}=U^{\odot n}$.

\proclaim{Proposition 3.4}
Let $V$ be  a unitary operator with a continuous spectrum.
Let $\sigma$ denote a measure of maximal spectral type of $V$.
The following holds.
\roster
\item"\rom{(i)}"
$\Cal M(V^{\otimes n})=n!\Cal M(V^{\odot n})$ (see, e.g., \cite{Ka}, \cite{KaLe}, \cite{Ry5}).
\item"\rom{(ii)}"
If $V^{\odot n}$ has a simple spectrum then
$\Cal M(V^{\otimes n/\Gamma})=\{n!/\#\Gamma\}$
for each subgroup $\Gamma\subset\goth S_n$.
In particular,
$\Cal M(V^{\otimes n})=\{n!\}$ and $\Cal M(V^{\odot (n-1)}\otimes V)=\{n\}$ (see, e.g., \cite{Ka}, \cite{Ag7}, \cite{DaRy1}).
\item"\rom{(iii)}"
If $V^{\odot n}$ has a simple spectrum then
 $V^{\odot k}$ has a simple spectrum for
 each $1\le k<n$ and $\sigma^{*i}\perp\sigma^{*j}$
 whenever $i\ne j$ and $i+j\le n$ \cite{Le4, Chapter 3}.
\item"\rom{(iv)}"
If $V$ has a simple spectrum and $\text{\rom{WCP}}(V)\ni a_iI+b_iV$ with $a_i,b_i\in\Bbb C\setminus\{0\}$ and
$\#\{a_1/b_1,\dots,a_n/b_n\}=n$
then $V^{\odot n}$
has a simple spectrum and $\sigma^{*i}\perp \sigma^{*j}$ whenever
$1\le i<j\le n$ (see, e.g., \cite{KaLe}, \cite{DaRy1}, \cite{Ag7}, \cite{Ry4}).
\endroster
\endproclaim

 Let $T$ be a transformation of $(X,\mu)$.
 Then $\goth S_n$ acts on the product space $(X,\mu)^{ n}$ by permutations of coordinates.
Moreover, this action commutes with the Cartesian product  $T^{\times n}$.
Hence for each subgroup $\Gamma\subset\goth S_n$, the $\sigma$-algebra of $\Gamma$-fixed  subsets in $X^n$ is invariant under $T^{\times n}$, i.e. it a {\it factor} of $T^{\times n}$.
We denote the restriction of $T^{\times n}$ to this $\sigma$-algebra by
$T^{\times n}/\Gamma$.

If $k<n$ then we think of $\goth S_k$ as the subgroup
$\{\tau\in\goth S_n\mid\tau(i)=i\text{ for all }i>k\}$ of $\goth S_n$.

\proclaim{Theorem 3.5 \cite{Ag7}}
Let $T$ be a weakly mixing transformation
on $(X,\goth B,\mu)$ with a simple spectrum.
If $a_iI+b_iU_T\in\text{\rom{WCP}}(U_T)$ for some $a_i,b_i\in\Bbb C\setminus\{0\}$, $i=1,\dots,n$, and
$\#\{a_1/b_1,\dots,a_n/b_n\}=n$
then  $\Cal M(T^{\times n})=\{n,n(n-1),\dots,n!\}$.
More generally, given a subgroup $\Gamma\subset\goth S_n$,
$$
\Cal M(T^{\times n}/\Gamma)=\{\#(\Gamma\backslash\goth S_n/\goth S_{k})\mid k=1,\dots, n-1\}.
$$
 In particular, $\Cal M(T^{\times n}/\goth S_{n-1})=\{2,3,\dots,n\}$.
\endproclaim

\demo{Proof}
We let $\Cal H:=L^2_0(X,\mu)$.
Since $L^2(X,\mu)=\Bbb C\oplus{\Cal H}$, it follows that
$$
\align
L^2(X,\mu)^{\otimes n}&=
\Bbb C\oplus \bigoplus_1^n{\Cal H}\oplus\bigoplus_1^{n(n-1)/2}{\Cal H}^{\otimes 2}\oplus\cdots\oplus{\Cal H}^{\otimes n},
\tag3-1\\
U_{T^{\times n}}&=\bigoplus_1^nU_T \oplus\bigoplus_1^{n(n-1)/2}U_T^{\otimes 2}\oplus\cdots\oplus U_T^{\otimes n}.
\tag3-2
\endalign
$$
It now follows from Proposition~3.4(i), (iv) that
$\Cal M(T^{\times n})=\{n,n(n-1),\dots,n!\}$.

The second claim is proved in a similar way.
\qed
\enddemo

\example{Example 3.6 \cite{Ag4}} For  the following non-mixing rank-one transformations $T$ and each $n>0$, the $n$-th symmetric power $T^{\odot n}$ of $T$ has a simple spectrum  (equivalently, the unitary operator $\exp U_T:=\bigoplus_{n=0}^\infty U_T^{\odot n}$ has a simple spectrum) and hence
$\Cal M(T^{\times n})=\{n, n(n-1),\dots,n!\}$:
 \roster
 \item"(i)"
 $T$ is the  Chacon  transformation
 associated with the sequence $(r_n,s_n)_{n=1}^\infty$, where
 $r_n\equiv 3$ and $s_n(0)=s_n(2)=0$ and $s_n(1)=1$.\footnote{This result was preceded by \cite{PrRy}, where it was shown that the convolutions of the maximal spectral type of Chacon transformation are pairwise disjoint.}
 \item"(ii)"   $T$ is del Junco-Rudolph map from \cite{dJRu}.
 It is associated with $(r_n,s_n)_{n=1}^\infty$, where $r_n=2^{n+1}$, $s_n(2^n)=1$ and $s_n(i)=0$ if $i\ne 2^n$.
 \endroster
  \endexample

The next natural problem is to investigate the spectral multiplicities of
the products $T_1^{n_1}/\Gamma_1\times\cdots\times T_k^{n_k}/\Gamma_k$, where $\Gamma_i$ is a subgroup in $\goth S_{n_i}$, $i=1,\dots,k$.
In this connection  an important property of unitary operators
 was  introduced in \cite{Ry5}.

\definition{Definition 3.7} Let $U$ and $V$ be unitary operators and let $\sigma_U$ and $\sigma_V$ denote  measures  of maximal spectral type of $U$ and $V$.
We say that $U$ and $V$ are {\it strongly disjoint} if the  map
$
(\Bbb T\times\Bbb T,\sigma_U\times\sigma_V)\ni (z_1,z_2)\mapsto z_1z_2\in (\Bbb T,\sigma_U*\sigma_V)
$
is one-to-one mod $0$.
\enddefinition

If $U$ and $V$ are strongly disjoint with continuous spectrum then their maximal spectral types are orthogonal.
The converse is not true.
For instance, if $U$ has Lebesgue spectrum then $U$ is not strongly disjoint from any $V$ with continuous spectrum.

\proclaim{Proposition 3.8} Let $U$ and $V$ be two unitary operators with  simple spectrum.
\roster
\item"\rom{(i)}" $U$ and $V$ are strongly disjoint if and only if
$U\otimes V$ has a simple spectrum.
\item"\rom{(ii)}" \cite{Ry5}
If $\text{\rom{WCP}}(U\otimes V)\ni aU\otimes I$ for some $0\ne a\in\Bbb C$ then $U$ and $V$ are strongly disjoint.
\item"\rom{(iii)}" \cite{KaLe}
If  $\text{\rom{WCP}}(U\otimes V)\supset\{ a(I+U)\otimes (I+V),
 a(I+U)\otimes (I+bV)\}$ for some $0\ne a,b\in\Bbb C$ and $b\ne 1$ then $U$ and $V$ are strongly disjoint.
\endroster
\endproclaim

 It follows from Proposition~3.8(i) that   if $U$ and $V$ are arbitrary strongly disjoint unitary operators then $\Cal M(U\otimes V)=\Cal M( U)\cdot\Cal M(V)$.
As in \cite{Ry5}, given subsets $E,F\subset\Bbb N$,
we let $E\diamond F:=E\cup F\cup E\cdot F$.

\proclaim{Corollary 3.9}
Let $T$ and $S$ be two weakly mixing transformations.
If $T$ and $S$ are strongly disjoint and  the  maximal spectral types of $U_T\otimes U_S$ and $U_T\oplus U_S$ are orthogonal
(for instance, if  there are non-zero $a_1,a_2,a_3$ such that \rom{WCP}$(U_T\oplus U_S)\ni\{0\oplus a_1I, a_2I\oplus 0, a_3U_T\oplus I\}$) then
then $\Cal M(T\times S)=\Cal M(T)\diamond\Cal M(S)$.
\endproclaim

\proclaim{Theorem 3.10 \cite{Ry5}}
There exist weakly mixing transformations  $T_1,T_2,\dots$ such that
the unitary operator $\exp(U_{T_1})\otimes\exp(U_{T_2})\otimes\cdots$ has a simple spectrum.
Hence
$$
\Cal M(T_1^{\times n_1}/\Gamma_1\times\cdots\times T_k^{\times n_k}/\Gamma_k)=
\Cal M(T_1^{\times n_1}/\Gamma_1)\diamond\cdots\diamond\Cal M( T_k^{\times n_k}/\Gamma_k).
$$
for each  finite sequence of positive integers $n_1,\dots,n_k$ and subgroups  $\Gamma_1\subset\goth S_{n_1}$,\dots, $\Gamma_k\subset\goth S_{n_k}$.
A similar assertion holds also for  infinite sequences
 $n_1,n_2,\dots$ and $\Gamma_1,\Gamma_2,\dots$.
\endproclaim

  It follows, in particular, that for each $m>0$, the set $\{m,m+1,m+2,\dots\}$ is realizable as the set of spectral multiplicities of the infinite product transformation $T_1^{\times m}\times T_2^{\times(m+1)}\times T_3^{\times(m+2)}\times\cdots$.

\head 4. Realization of subsets containing $2$
\endhead

In this section we explain how  to combine
 the technique of isometric extensions
 (Section~2)
with the technique of Cartesian products  (Section~3)
to realize the subsets containing $2$.

\proclaim{Theorem 4.1}
Given a subset $E\subset\Bbb N$, there is a weakly mixing transformation
$R$ with $\Cal M(R)=E\cup\{2\}$.
\endproclaim
\demo{Proof}
For finite $E$, this theorem was proved in  \cite{KaLe}.
 For arbitrary $E$, it was proved in \cite{Da3}  and  \cite{Da4}.
Let $(X,\mu,T),D,K,\Cal K,Y$ and $\tau=(\psi,\phi)$  denote the same objects as in the proof  of Theorem~2.1.
We recall that $T$ is a rank-one transformation associated with a sequence $(r_n,s_n)_{n=1}^\infty$ and $\tau$ is a product-type cocycle of $T$
associated with a sequence $(d_n,b_n)_{n=1}^\infty$.
Some conditions were imposed on the sequence $(r_n,s_n,d_n,b_n)_{n=1}^\infty$ to obtain
$\Cal M(T_{Y\backslash\tau})=E\cup\{1\}$.
We will  add some more conditions to construct the desired transformation.
\proclaim{Claim 4.2}
Let $k\in\Cal K$ and $n_l\to\infty$.
If $r_{n_l}\to\infty$, $s_{n_l}(j)=0$ for  $0\le j<r_{n_l}/2$ and
$s_{n_l}(j)=1$ for  $r_{n_l}/2\le j<r_{n_l}$,
$d_{n_l}\equiv 0$,
$b_{n_l}(j)=0$ for  $r_{n_l}/2\le j<r_n$
and for each $q\in \Cal O_*(k)$,
$$
\frac{\#\{j\mid 0\le j<r_{n_l}/2,b_{n_l}(j)-b_{n_l}(j+1)=q\}}{r_{n_l}}\to\frac 1{2\#\Cal O_*(k)}
$$
then
$$
U_{T_\tau,\chi}^{h_{n_l}}\to 0.5(l_\chi(k)I+U_{T_\tau,\chi}^*)\text{ for each }\chi\in\widehat K\text{  as }l\to\infty.\tag4-1
$$
\endproclaim

This claim is proved in a similar way as Claim 2.8.
Thus from now on we will assume that  for each $k\in\Cal K$ there is a sequence $n_l\to\infty$ such that \thetag{4-1} holds.

Let $R:=T_{Y\backslash\tau}\times T_\psi$.
Then
$$
U_{R}=U_{T_\psi\times T_\psi}\oplus\bigoplus_{1\ne\chi\in \widehat{K/Z}} U_{T_\tau,\chi}\otimes (U_{T_\psi}\oplus P_0),\tag4-2
$$
where $P_0$ is the orthogonal projection to the subspace of constants in $L^2(X\times D,\mu\times\lambda_D)$.

{\bf a)} Since $U_{T_\psi}$ has a simple continuous spectrum and WCP$(U_{T_\psi})\ni 0.5(I+U_{T_\psi}^*)$, it follows from Theorem~3.3 that  $\Cal M(T_\psi\times T_\psi)=\{2\}$.

{\bf b)} For each $\chi\in\widehat{Y\backslash K}$, the unitary operator
$U_{T_\tau,\chi}$ has a simple spectrum.
If $\chi\ne 1$, there is $k\in\Cal K$ with $\chi(k)\ne 1$.
It follows from \thetag{4-1} that WCP$(U_{T_\tau,\chi}\otimes (U_{T_\psi}\oplus P_0))$ contains
$0.25(l_\chi(k)I+U_{T_\tau,\chi}^*)\otimes((I+U_{T_\psi}^*)\oplus P_0)$
and $0.25(I+U_{T_\tau,\chi}^*)\otimes((I+U_{T_\psi}^*)\oplus P_0)$.
Therefore by Proposition~3.8(iii),
$U_{T_\tau,\chi}^*$ and $U_{T_\psi}^*\oplus P_0$ are strongly disjoint
\footnote{Slightly abusing notation we use that $I\oplus P_0=I$.}.
Hence the unitary operator $U_{T_\tau,\chi}^*\otimes(U_{T_\psi}^*\oplus P_0)$ has a simple spectrum by Proposition~3.8(i).

{\bf c)}
If $\chi,\chi'\in\widehat{Y\backslash K}$ and $\chi\approx\chi'$ then
 $U_{T_\tau,\chi}$ and $U_{T_\tau,\chi'}$
 are unitarily equivalent.
Hence
$U_{T_\tau,\chi}\otimes (U_{T_\psi}\oplus P_0)$ and $U_{T_\tau,\chi'}\otimes (U_{T_\psi}\oplus P_0)$ are also unitarily equivalent.

{\bf d)}
If $\chi,\chi'\in\widehat{Y\backslash K}$ but $\chi\not\approx\chi'$
then there is $k\in\Cal K$ such that $l_\chi(k)\ne l_{\chi'}(k)$.
Moreover, by Claim~2.8, there is a sequence $h_{n_l}\to\infty$ such that
$U_{T_\psi}^{h_{n_l}}\to I$, $U_{T_\tau,\chi}^{h_{n_l}}\to l_\chi(k)I$,  and $U_{T_\tau,\chi}^{h_{n_l}}\to l_{\chi'}(k)I$.
It follows from Lemma~2.2 that the maximal spectral types of
$U_{T_\tau,\chi}\otimes (U_{T_\psi}\oplus P_0)$ and $U_{T_\tau,\chi'}\otimes (U_{T_\psi}\oplus P_0)$ are orthogonal.

We now deduce from {\bf a), b), c), d)} and \thetag{4-2} that
$\Cal M(R)=E\cup\{2\}$.
\qed
\enddemo

\remark{Remark \rom{4.3}} We encounter here with an interesting phenomenon.
The transformation $R$ can be represented an an isometric extension
$(T_\psi\times T_\psi)_{\upsilon}$ of the product
$T_\psi\times T_\psi$ which has a homogeneous spectrum of multiplicity $2$, where $\upsilon:=(Y\backslash\widetilde\phi)\otimes 1$.
However the unitary operators $U_{(T_\psi\times T_\psi)_\upsilon,\chi}$
related to the extension (as in (2-1)) have simple spectrum if $\chi\ne 1$.
This  property of the extension (appeared first in \cite{KaLe}) was possible to realize due to the fact that the cocycle $\upsilon$ depends only on the ``first coordinate'', i.e. $\upsilon$ is, in fact, a cocycle of the first marginal factor $T_\psi$ of $T_\psi\times T_\psi$.
\endremark

The following theorem is a  generalization of Theorem~2.10.
It is proved in a similar way.

\proclaim{Theorem 4.4} The subset
of pairs
$(T,\tau)\in \text{\rom{Aut}}(X,\mu)\times\Cal A(X,Z)$ such that
$T_{Y\backslash\tau}$ is weakly mixing,  $\Cal M (T_{Y\backslash\tau})=E\cup\{1\},\
\Cal M (T_{Y\backslash\tau}\times T_{\psi})=E\cup\{2\}
$
is residual in $\text{\rom{Aut}}(X,\mu)\times\Cal A(X,Z)$.
\endproclaim

\remark{Remark \rom{4.5 }}
As was noted in \cite{KaLe} and  \cite{Da3},
a more general class of subsets can be realized as spectral multiplicities of transformations
$
T_{Y\backslash\tau}\times (T_\psi)^{\times k}/\Gamma,
$
where $\Gamma$ is a subgroup of $\goth S_k$, $k>1$.
However for that we need to modify the cocycle $\tau$ in such a way
that the transformation
$(T_\psi)^{\odot(k+1)}$ has a simple spectrum.
This can be done by manufacturing several additional weak limits for
$U_{T_\psi}$ and applying
Theorem~3.5.
 In particular, setting $k=2$ and $\Gamma$ trivial we obtain that {\it every}
subset $\{3,6\}\cup F$, where $F$ is an arbitrary subset of even numbers, is realizable.
\endremark

\head 5. Rokhlin problem and auxiliary group actions
\endhead

It was shown in  Theorem~3.3  that  the Cartesian square of a weakly mixing transformation $T$ on $(X,\goth B,\mu)$ with a simple spectrum has a homogeneous spectrum of multiplicity $2$ (under some condition).
This answers (Pr3) in the particular case $n=2$.
The proof of Theorem~3.3 (or, more precisely, Lemma~3.1, which is the main ingredient of that proof) uses implicitly the fact that the transformation $T\times T$ embeds into the action $W=(W_g)_{g\in G_2}$ of the non-Abelian group $G_2:=\Bbb Z^2\rtimes (\Bbb Z/2\Bbb Z)$ on $(X\times X,\mu\times\mu)$ generated by
$$
W_{(n,m,0)}(x,y)=(T^nx,T^my) \quad\text{ and }\quad W_{(0,0,1)}(x,y)=(y,x),\tag5-1
$$
where $n,m\in\Bbb Z$.
Indeed, $T\times T=W_{(1,1,0)}$.
This hints that in order to obtain transformations with homogeneous spectrum of multiplicity $n>2$, it it natural to consider actions of the meta-Abelian group $G_n:=\Bbb Z^{n}\rtimes\Bbb Z/n\Bbb Z$, where the cyclic group $\Bbb Z/n\Bbb Z$ acts on $\Bbb Z^{n}$ by cyclic permutations.
Of course, given a transformation $T$ of $(X,\mu)$, we can form a $G_n$-action on the product space $(X,\mu)^{\times n}$ generated by two transformations
$$
\aligned
(x_1,x_2,\dots,x_n)&\mapsto(Tx_1,x_2,\dots,x_n) \quad \text{and}\\
(x_1,x_2,\dots,x_n)&\mapsto(x_2,x_3,\dots,x_1).
\endaligned
\tag5-2
$$
However as follows from Theorem~3.5 the transformation $T^{\times n}$
(which is included into the action of $G_n$ as the $n$-th power of the composition of the two transformations defined by \thetag{5-2})
does not have homogeneous spectrum.
The main reason for that is
the
presence of the one-dimensional invariant subspace of constants in $L^2(X,\mu)$.
Because of it we obtain {\it long sums} \thetag{3-1} and  \thetag{3-2}, where every summand contributes  non-trivially  into the set of spectral multiplicities of $T^{\times n}$.
However  there exist  other actions of $G_n$ that help to solve the  Rokhlin problem.
We call them {\it suitable auxiliary} actions.
It is not an easy task to  construct them in an explicit  way (the only known effective cutting-and-stacking construction of such actions is given in  \cite{Da1}).
It is much easier to prove their existence via Baire  category methods.
As appears, a generic $G_n$-action is suitable.

We need some notation.
Each element of $G_n$ is a sequence $(k_1,\dots,k_{n+1})$, where $k_1,\dots,k_n\in\Bbb Z$ and  $k_{n+1}\in\Bbb Z/n\Bbb Z$.
Let $e_1:=(1,0,\dots,0),\dots, e_n:=(0,\dots,1,0)$, $e_*:=(1,\dots,1,0)$ and  $e_0:=(0,\dots,0,1)$.
Then $e_*=(e_1e_0)^n$.
By $H_n$ we denote the subgroup of $G_n$ generated by $e_1,\dots,e_n$.
Let $A$ denote the automorphism of $H_n$ generated by the conjugation by $e_{0}$.

\proclaim{Proposition 5.1}
Let $U=(U_g)_{g\in G_n}$ be a unitary representation of $G_n$ in a Hilbert space $\Cal H$.
If the operator $U_{e_{1}e_0}$ has a simple continuous spectrum and $U_{e_je_1^{-1}}$
has a continuous spectrum for each $2\le j<n$ then $U_{e_{*}}$ has a homogeneous continuous spectrum of multiplicity $n$.
\endproclaim
\demo{Proof}
By the spectral theorem for $(U_h)_{h\in H_n}$,
$$
\Cal H=\int^\oplus_{\widehat {H_n}}\Cal H_w\,d\sigma(w)
 \text{ \ and \  }U_h=\int^\oplus_{\widehat {H_n}}\langle h,w\rangle I_w\,d\sigma(w)\tag5-3
$$
for all $h\in H_n$, where $I_w$ stands for the identity operator in the Hilbert space $\Cal H_w$.
The inclusion $\Bbb Z\ni m\mapsto me_{*}\in H_n$
generates the canonical  projection $\widehat {H_n}\to\Bbb T$.
Let $\sigma=\int_\Bbb T\sigma_z\,d\widetilde\sigma(z)$
be the disintegration of $\sigma$ with respect to this projection.
Then we derive from \thetag{5-3} that
$$
\Cal H=\int^\oplus_{\Bbb T}\Cal H_z'\,d\widetilde\sigma(z)
 \ \text{ and }\ U_{e_{*}}=\int^\oplus_{\widehat {H_n}}z I_z\,d\widetilde\sigma(w),
$$
where $\Cal H_z'=\int_{\widehat {H_n}}^{\oplus}\Cal H_w\,d\sigma_z(w)$.
Denote by $\widehat A$ the automorphism of $\widehat {H_n}$ which is dual to $A$.
Since the unitary representation $(U_{Ah})_{h\in H_n}$
is unitarily equivalent to  $(U_{h})_{h\in H_n}$,
it follows that dim$\,\Cal H_w=\text{dim}\,\Cal H_{\widehat Aw}$
 for a.a. $w$ and $\widetilde\sigma\circ \widehat A$ is equivalent to $\widetilde \sigma$.
Without loss of generality we may assume that
$\widetilde\sigma\circ \widehat A=\widetilde \sigma$.
 This implies that $\sigma_z\circ \widehat A=\sigma_z$ for a.a. $z$.
Hence dim$\,\Cal H_z'=\text{dim}\,\Cal H_{\widehat Az}'$ for a.a. $z$.
We now claim that a.e.  $\widehat A$-orbit consists of $n$ points.
Indeed, this follows from the fact that  $U_{e_{j}e_1^{-1}}$ has no non-trivial fixed vectors because the latter implies
$\sigma(\{z=(z_1,\dots,z_n)\in\widehat{H_n}\mid z_1=z_j\})=0$  for  each $j>1$.
Therefore if $\sigma_z$ has an atom at some point $w$ then  it  has also atoms at $n-1$ points $\widehat A^jw$, $j=2,\dots,n$.
Hence $\Cal M(U_{e_{*}})\subset\{n,2n,\dots\}\cup\{\infty\}$.
Since $U_{e_1e_0}$ has a simple spectrum and
$U_{e_1e_0}^n=U_{e_{*}}$, it follows that $\Cal M(U_{e_{*}})\subset\{1,\dots,n\}$.
Hence $\Cal M(U_{e_{*}})=\{n\}$.
\qed
\enddemo

Denote by $\Cal A_{G_n}$ the space of all measure preserving $G_n$-actions on $(X,\mu)$.
 We consider $\Cal A_{G_n}$ as a subset of the infinite product space $\text{Aut}(X,\mu)^{G_n}$ endowed with the product (Polish) topology.
Then $\Cal A_{G_n}$ is closed and hence Polish in the induced topology.
There is a natural continuous action of  Aut$(X,\mu)$  on $\Cal A_{G_n}$
by conjugation:
$$
(R\cdot T)_g:=RT_gR^{-1}.
$$

\proclaim{Theorem 5.2} Fix $n>1$, a sequence $k_l\to\infty$ and  a polynomial $p$ with non-negative coefficients and $p(1)\le1$.
Then the subset of $G_n$-actions $T=(T_g)_{g\in G_n}$ such that
the transformation $T_{e_{*}}$ is weakly mixing, $\Cal M(T_{e_{*}})=\{n\}$ and $(U_{T_{e_{*}}})^{-m_l}\to p(U_{T_{e_{*}}})$ along a subsequence $(m_l)_{l=1}^\infty$ of $(k_l)_{l=1}^\infty$
is residual in $\Cal A_{G_n}$.
\endproclaim

\demo{Proof}
The subsets
$$
\align
\Cal W_{G_n}&:=\{T\in \Cal A_{G_n}\mid T_g\text{ is weakly mixing  for each $g\in G_n$ of infinite order}\},\\
\Cal S_{G_n}&:=\{T\in \Cal A_{G_n}\mid T_{e_1e_0}\text{ has a simple spectrum}\},\\
\Cal L_{G_n}&:=\{T\in \Cal A_{G_n}\mid (U_{T_{e_{*}}})^{-m_l}\to p(U_{T_{e_{*}}})
\text{ for some $(m_l)_{l}\subset(k_l)_{l}$}
\}
\endalign
$$
are $G_\delta$ in $\Cal A_{G_n}$ because the subset of weakly mixing transformations, the subset of transformations with  simple spectrum and the subset
$$
\Cal J:=\{T\in\text{Aut}(X,\mu)\mid
(U_{T_{e_{*}}})^{-m_l}\to p(U_{T_{e_{*}}})
\text{ for some $(m_l)_{l}\subset(k_l)_{l}$}
\}
$$
 are all $G_\delta$ in Aut$(X,\mu)$.
The three subsets $\Cal W_{G_n}$, $\Cal S_{G_n}$ and $\Cal L_{G_n}$ are invariant under  conjugation.
We recall that the Aut$(X,\mu)$-orbit of every free $G_n$-action is dense in $\Cal A_{G_n}$ \cite{FoWi}.
 Since each Bernoulli $G_n$-action is free and belongs to $\Cal W_{G_n}$, it follows that  $\Cal W_{G_n}$ is a dense $G_\delta$ in $\Cal A_{G_n}$.

Let us show that $\Cal S_{G_n}$ contains a free $G_n$-action.
We define  an action $T$ on the product space $(Y,\nu):=(\Bbb T,\lambda_\Bbb T)\times(\Bbb Z/n\Bbb Z,\lambda_{\Bbb Z/n\Bbb Z})$ by
setting
$T_{e_1}(x,i):=(S_ix, i)$ and $T_{e_0}(x,i):=(x,i\oplus 1)$,
where $\oplus$ means addition mod $n$ and
$S_1,\dots,S_n$ are rotations by rationally independent irrationals.
Since $e_0$ and $e_1$ generate $G_n$, the $G_n$-action $T$ is well defined.
This action is free.
The transformation $T_{e_0e_1}$ is ergodic and has a simple spectrum.
Hence $T\in \Cal S_{G_n}$.
We thus obtain that  $\Cal S_{G_n}$ is a dense $G_\delta$ in $\Cal A_{G_n}$.

Take an aperiodic transformation $S\in\Cal J$ and place it into \thetag{5-2} instead of $T$.
Then  \thetag{5-2} defines a free $G_n$-action which belongs to $\Cal L_{G_n}$.
Hence $\Cal L_{G_n}$ is a dense $G_\delta$ in $\Cal A_{G_n}$.

The intersection
$\Cal S_{G_n}\cap\Cal W_{G_n}\cap\Cal L_{G_n}$
is also  a dense $G_\delta$ in $\Cal A_{G_n}$.
It remains to note that if $R\in\Cal S_{G_n}\cap\Cal W_{G_n}$
then $\Cal M(R_{e_{*}})=\{n\}$ by Proposition~5.1.
\qed
\enddemo

\remark{Remark \rom{5.3}} It is possible also to take  other {\it auxiliary} meta-Abelian groups instead of $G_n$ to answer the Rokhlin question on homogeneous spectrum.
For instance, in the original paper \cite{Ag6} Ageev considered the groups
$\Bbb Z^{n-1}\rtimes_{A'}\Bbb Z$, where $A'$ is a group automorphism of $\Bbb Z^{n-1}$
defined by $A'e_1:=e_2,\dots,A'e_{n-2}:=e_{n-1}$ and
$A'e_{n-1}:=-e_1-\cdots-e_{n-1}$.
The present author used the groups $G_n':=\Bbb Z^{n}\rtimes_{A}\Bbb Z$ in \cite{Da1}, where $A$ is the {\it cyclic}  automorphism  of $\Bbb Z^{n}$ considered
in this section.
One more version for $G_n$ is suggested in \cite{DaSo}.
In these notes we choose $G_n$  from \cite{Ry5}  for the only reason:
 in this case Proposition~5.1 looks as  a  natural extension of Lemma~3.1
 and hence the method of suitable auxiliary actions is a natural development of the method of Cartesian powers.
\endremark

The following can be deduced from Theorem~5.2 and Corollary~3.9.

\proclaim{Theorem 5.4 \cite{Ry5}}
\roster
\item"\rom{(i)}"
Let $S$ be a rigid weakly weakly mixing
transformation.
Then for each $p>1$, there is a weakly mixing transformation $T$ such that $\Cal M(S\times T)=\Cal M(S)\diamond\{p\}$.
\item"\rom{(ii)}"
Given a finite set $p_1,\dots,p_l\in\Bbb N$, there is a weakly mixing transformation $T$ with $\Cal M(T)=\{p_1\}\diamond\cdots\diamond\{p_l\}$.
A similar statement holds also  for an infinite set $\{p_1,p_2,\dots\}$.
\endroster
\endproclaim

In particular, the sets $\{p,q,pq\}$,  $\{p,q,r,pq,pr,qr, pqr\}$, \dots are realizable for arbitrary positive integers $p,q,r,\dots$.

\proclaim{Corollary 5.5 \cite{Ry5}}
Each multiplicative semigroup of positive integers  is realizable.
Each additive semigroup of positive integers is realizable.
\endproclaim

We note that every additive subsemigroup of $\Bbb N$ is a multiplicative subsemigroup of $\Bbb N$.

\comment
The next assertion refines Theorem~7.10.
It can be deduced  from Propositions~7.4(iii), 7.8(iii) and Theorem~9.2.

\proclaim{Corollary 9.7?} For each sequence of positive integers $m_1,m_2,\dots$ there exist weakly mixing transformations $T_1,T_2,\dots$
such that $\Cal M(T_1^{m_1})=\{m_1\}$,
$\Cal M(T_1^{m_2})=\{m_2\}$,\dots
and the operator $\exp(U_{T_1}\oplus U_{T_2}\oplus\cdots)$
has a simple spectrum.
\endproclaim
\endcomment

The following remark was communicated to the author by V.~Ryzhikov.

\remark{Remark \rom{5.6}}
It follows from Proposition~3.4(iii) that if  $\exp(U_T)$ has a simple spectrum then the convolutions
of $\sigma_T$ are pairwise disjoint.
The converse is not true.
Indeed, by Theorem~5.2,  there is a weakly mixing transformation $T$ such that $\Cal M(T)=\{n\}$ and WCP$(U_T)\ni\{\alpha I+(1-\alpha)U_T\mid 0<\alpha<1\}$.
 Proposition~3.4(iv) yields that $\Cal M(U_{T}\odot U_{T})=\{n^2\}$.
On the other hand, again by Proposition~3.4(iv), $\sigma_{T}^{*i}\perp\sigma_{T}^{*j}$ for all $i\ne j$.
\endremark

 In the following theorem the technique of auxiliary group actions is combined with the technique of isometric extensions to  obtain spectral realizations of a new class of subsets in $\Bbb N$.

\proclaim{Theorem 5.7 \cite{Da1}}
Given  $n>0$ and a subset $E\subset\Bbb N$,
there is a weakly mixing  transformation $T$
with $\Cal M(T)=n\cdot (E\cup\{1\})$.
\endproclaim
\demo{Idea of the proof}
We start with a explicit construction of a weakly mixing $R$ with $\Cal M(R)=\{n\}$ from \cite{Da1}.
For that we produce via cutting-and-stacking an auxiliary  group action of $G_n'$ (defined in Remark~5.3).
Then $R=S^n$, where $S$ is a transformation of rank one.
By Lemma~A.3, there exist a compact Abelian group $K$, a closed subgroup $Y$ of $K$ and an automorphism $v$ of $K$ such that
$E=L(\widehat K,\widehat v,\widehat{Y\backslash K})$.
While doing cutting and stacking we construct also an ergodic cocycle $\phi$ of $S$ with values in $K$ such that
\roster
\item"$\bullet$" the operator
$U_{S_{\phi},\chi}$ has a simple spectrum for each $\chi\in\widehat K$,
\item"$\bullet$"
 maximal spectral types of the unitary operators $(U_{S_{\phi},\chi})^n$ and $(U_{S_{\phi},\chi'})^n$ are  orthogonal if $\chi'$ does not belong to the $\widehat v$-orbit of $\chi$.
\item"$\bullet$"
$U_{S_{\phi},\chi}$ and
$U_{S_{\phi},\chi'}$ are unitarily equivalent
$\chi'$  belongs to the $\widehat v$-orbit of $\chi$.
\endroster
These properties imply that $\Cal M(S_{Y\backslash\phi})=E\cup\{1\}$.
The construction retains enough  freedom to satisfy some extra conditions
(listed in Proposition~5.1)
which guarantee that  $\Cal M((U_{S_{\phi},\chi})^n)=\{n\}$ for each $\chi\in\widehat K$.
Since $(U_{S_{\phi},\chi})^n=U_{R_{\phi^{(n)}},\chi}$, it follows that
 $\Cal M(R_{Y\backslash\phi^{(n)}})=n\cdot(E\cup\{1\})$.
 \qed
\enddemo

\head 6. Mixing realizations
\endhead

\subhead Almost staircase systems
\endsubhead
These dynamical  systems play an important role in constructing mixing transformations with non-simple spectrum.

\definition{Definition 6.1} Let $T$ be a rank-one transformation associated with a sequence $(r_n,s_n)_{n=1}^\infty$.
If $s_n(i)=i$ for all $0\le i<r_n$ then $T$ is called a {\it staircase} transformation.
If there is a sequence $\delta_n\to 0$ such that $s_n(i)=i$ for all $\delta_nr_n\le i<r_n$ then
$T$ is called an {\it almost staircase} transformation.
\enddefinition

Smorodinsky conjectured that the classical\footnote{ The classical staircase corresponds to the case $r_n=n$.} staircase construction
  is mixing.
Adams \cite{Ad} proved the Smorodinsky conjecture: if $T$ is a finite measure preserving staircase transformation such that $r_n\to\infty$ and $r_n^2/h_n\to 0$
 then $T$ is mixing.
 We recall that $h_n$ stands for the hight of the $n$-th tower.
The condition $r_n^2/h_n\to 0$ is called {\it the restricted growth condition}.
Adams asked if it is possible to remove it from the statement of his theorem?
The affirmative answer was announced by Ryzhikov in 2000.
A detailed proof of this fact appeared recently in \cite{CrSi}.
It follows from Adams theorem and its generalization in \cite{CrSi}
 that the finite measure preserving almost staircase transformations are also mixing (with or without the restricted growth condition).

\subhead To force mixing
\endsubhead
Suppose that we are given a subset $E\subset \Bbb N$  to  construct a mixing transformation $T$ with $\Cal M(T)=E$.
For that we construct $T$ as a {\it limit} of a certain sequence of  non-mixing transformations $T_n$ with $\Cal M(T_n)=E$.
The transformations $T_n$ consist of a {\it rigid part}  and a {\it mixing part}.
The mixing part {\it occupies more and  more space} when $n\to\infty$.
Therefore the limit $T$ of this sequence is a mixing transformation.
The rigid part is needed to produce a ``rich'' semigroup
WCP$(U_{T_n})$ which, in turn, is used to control the  spectral multiplicities $\Cal M(T_n)$.
 The control is implemented in the same way as in  Sections~2--4.
 Then an approximation technique (based on Lemma~2.3) is used to retain the property $\Cal M(T_n)=E$ in the limit, i.e. to obtain $\Cal M(T)=E$.
 We note that the technique of weak limits can not be used straightforwardly to control $\Cal M(T)$ because $T$ is mixing and hence WCP$(U_T)=\{U_T^j\mid j\in\Bbb Z\}\cup\{0\}$.

We illustrate the method to force mixing by the following theorem.

\proclaim{Theorem 6.2 (\cite{Ag7}, \cite{Ry4})}
There is a mixing almost staircase transformation $T$ such that the unitary operator $\exp(U_T)$
has a simple spectrum.
\endproclaim

\demo{Idea of the proof}
We need a notation. Given a Hilbert space $\Cal H$, a vector $h\in\Cal H$, an operator $V$ in $\Cal H$ and an integer $K>0$, we denote by $\Cal L(V,h,K)$ the linear span of $\{V^jh\mid |j|\le K\}$.

 Fix $0<\delta<1$.
Let $T_\delta$ be a rank-one transformation associated with a sequence $(r_n,s_n)_{n=1}^\infty$ such that
 $r_n\to\infty$, $s_n(j)=j$  for $\delta r_n\le j<r_n$ and $n\ge 1$ and
 WCP$(U_{T_\delta})\supset\{\alpha I+\beta U_{T_\delta}\mid \alpha,\beta>0, \alpha+\beta<\delta\}$.
 The first tower of $T_\delta$  is not specified yet.
It follows from Proposition~3.4(iv) that the unitary operator $\exp(U_{T_\delta})$ has a simple spectrum.

Fix a sequence of positive reals $\delta_n\to 0$ and a sequence of integers
$N_n\to\infty$.
 Construct a sequence of rank-one transformations $T_{\delta_n}$ as above with an additional   {\it agreement} condition specifying their first towers:
 the first tower of $T_{\delta_{n+1}}$ coincides with the
 $N_n$-tower of $T_{\delta_n}$ for each $n>0$.
 We now define a new rank-one transformation $T$ via  a ``concatenation'' of the sequence $T_{\delta_n}$ in the following sense.
The first $N_1$ towers of $T$ are the first $N_1$ towers of $T_{\delta_1}$.
The next $N_2-1$ towers of $T$ are  the first $N_2-1$ towers of
 $T_{\delta_{2}}$.
 Then continue with  $N_3-1$ first towers of $T_{\delta_{3}}$ and so on.
 Due to the agreement condition this inductive cutting-and-stacking procedure   is well defined.
 The corresponding rank-one transformation $T$ is an almost staircase.
 It follows that $T$ is mixing.

It remains to explain how to select the sequence $(N_n)_{n=1}^\infty$  to obtain  simplicity of spectrum of  the operator $\exp (U_T)$.
 Let $X_n$ be the space of $T_{\delta_n}$ and let $\Cal H_n$ stand for the subspace of functions in $L^2_0(X_{n})$
 which are constant on every level of the $N_n$-tower and  vanish outside this tower.
Since  $\exp(U_{T_{\delta_n}})$ have a simple spectrum,
we can choose $N_n$  large so that  for each $s=1,\dots,n$, there is a vector $v_s\in\Cal H_n^{\odot s}$ and $K>0$ such that
\roster
\item"(i)"
$\Cal L(U_{T_{\delta_n}}^{\otimes s},v_s,K)\subset\Cal H_n^{\odot s}$
and
\item"(ii)"
dist$(v,\Cal L(U_{T_{\delta_n}}^{\otimes s},v_s,K))<\delta_n$
for each $v\in\Cal H_{n-1}^{\odot s}$ with $\|v\|=1$.
\endroster
We note that $\Cal H_n$ is also a subspace in $L^2_0(X_{n+1})$.
Moreover, $\Cal H_1\subset\Cal H_2\subset\cdots$ and the union $\bigcup_j\Cal H_j$ is dense in $L^2_0(X)$.
It follows from (i) that $\Cal L(U_{T_{\delta_n}}^{\otimes s},v_s,K)=\Cal L(U_T^{\otimes s},v_s,K)$.
Lemma~2.3 and (ii) now yield  that $U_T^{\odot s}$ has a simple spectrum for each $s$.
Hence $\exp( U_T)$ has a simple spectrum.
\qed
\enddemo

 Ryzhikov showed in \cite{Ry4} that there is also  a  mixing staircase transformation $T$ with non-monotone sequence $r_n\to\infty$ such that $\exp U_T$ has a simple spectrum.
He also constructed  a mixing $T$ such that the infinite product
$T\times T^2\times T^3\times\cdots$ has a simple spectrum \cite{Ry7}.
This transformation plays an important  role in Tikhonov's proof of Theorem~6.7 below (see \cite{Ti2}).

We now obtain  {\it mixing counterparts} of Theorems~3.3 and 3.5.

\proclaim{Corollary 6.3} There is a mixing transformation $T$ such that
\roster
\item"\rom{(i)}" \cite{Ry3}
$\Cal M(T\times T)=\{2\}$.
\item"\rom{(ii)}" \cite{Ag7}
More generally, given a subgroup $\Gamma\subset\goth S_n$,
$$
\Cal M(T^{\times n}/\Gamma)=\{\#(\Gamma\backslash\goth S_n/\goth S_{k})\mid k=1,\dots, n-1\}.
$$
 In particular, $\Cal M(T^{\times n}/\goth S_{n-1})=\{2,3,\dots,n\}$.
\endroster
\endproclaim

In a similar way we can force mixing for countably many pairwise strongly disjoint transformations  to obtain a {\it mixing} version of
Theorem~7.10.

\proclaim{Theorem 6.4 \cite{Ry5}}
There exist mixing transformations  $T_1,T_2,\dots$ such that
the unitary operator $\exp(U_{T_1})\otimes\exp(U_{T_2})\otimes\cdots$ has a simple spectrum.
Hence
$$
\Cal M(T_1^{\times n_1}/\Gamma_1\times\cdots\times T_k^{\times n_k}/\Gamma_k)=
\Cal M(T_1^{\times n_1}/\Gamma_1)\diamond\cdots\diamond\Cal M( T_k^{\times n_k}/\Gamma_k)
$$
for each  finite sequence of positive integers $n_1,\dots,n_k$ and subgroups  $\Gamma_1\subset\goth S_{n_1}$,\dots, $\Gamma_k\subset\goth S_{n_k}$.
A similar assertion holds also for  infinite sequences
 $n_1,n_2,\dots$ and $\Gamma_1,\Gamma_2,\dots$.
\endproclaim

Next, meta-Abelian extensions of rank-one maps from the proofs of Theorems~2.1  and 4.1 are also suited well to force mixing in them.

\proclaim{Theorem 6.5 \cite{Da4}} Let $E\subset\Bbb N$.
Then there are  mixing transformations $R_1$ and $R_2$ such that with $\Cal M(R_1)=E\cup\{1\}$ and
$\Cal M(R_2)=E\cup\{2\}$.

\endproclaim
\demo{Idea of the proof}
We discuss  only the construction of $R_1$.
As in the proof of Theorem~2.1, $R_1$ appears as a compact extension $T_{Y\backslash\tau}$, where $T$ is a rank-one transformation and  $\tau$ a cocycle of $T$ with values in $D\ltimes K$. The compact groups $D,K,Y$  are exactly same as in the proof of Theorem~2.1.
As in the proof of Theorem~6.2, $T_{Y\backslash\tau}$ appears a {\it limit} of a sequence of non-mixing weakly mixing
transformations $(T_n)_{Y\backslash\tau_n}$ such that $\Cal M((T_n)_{Y\backslash\tau_n})=E\cup\{1\}$.
 Every $T_n$ has a {\it rigid part} responsible for the weak limits in WCP$(U_{T_n})$ that jointly with the cocycles $\tau_n$ guarantee the desired spectral multiplicities of $(T_n)_{Y\backslash\tau_n}$
and  a {\it mixing}   (staircase) part.
 The latter  {\it occupies} $(1-\delta_n)$-part of the space where $T_n$ is defined and $\delta_n\to 0$.
Hence the limit $T$ of the sequence $T_n$ is an almost staircase and hence it is mixing.
The extension $T_{Y\backslash\tau}$ is weakly mixing by construction.
Hence it is mixing  \cite{Ru}.
 Lemma~2.3 is used to retain $E\cup\{1\}$ as the spectral multiplicity set in the {\it limit } of the sequence $(T_n)_{Y\backslash\tau_n}$.
Thus $\Cal M(R_1)=E\cup\{1\}$.
\qed
\enddemo

\remark{Remark \rom{6.6}}
\roster
\item"(i)"
In order to force mixing the authors of \cite{Ry3}, \cite{Ry4}, \cite{Ag7}, \cite{Da4} use Smorodin\-sky-Adams  staircase constructions.
 It is also possible to use stochastic Ornstein rank-one mixing constructions \cite{Or} in place of them.
\item"(ii)"
While {\it Abelian} compact extensions are convenient to construct weakly mixing realizations of subsets $E$ containing $1$ or $2$ (see the second proof of Theorem~2.1 and  \cite{Go--Li}, \cite{Da3}), it is unclear how to force mixing in them.
\item"(iii)"
All the aforementioned mixing realizations are mixing of all orders in view of \cite{Kal}, \cite{Ry1} and \cite{Ru}.
\endroster
\endremark

\subhead Generic approach
\endsubhead
The subset  Mix$(X,\mu)$ of mixing transformations  is meager in Aut$(X,\mu)$ endowed  with the weak topology.
 Therefore the weak topology is not suitable to apply the Baire category argument  in Mix$(X,\mu)$.
Tikhonov introduced another topology, say $\tau$, on Aut$(X,\mu)$ such that the subspace (Mix$(X,\mu),\tau)$ is Polish \cite{Ti1}.
Given  a subset $A\subset X$, we define a map
$f_A:\text{Aut}(X,\mu)\to l^\infty$ by setting
$f_A(T):=(\mu(T^nA\cap A))_{n\in\Bbb N}$.
Let $\tau$ be the weakest topology on $\text{Aut}(X,\mu)$ in which
all the maps $f_A$ are continuous.
It is assumed that the space $l^\infty$ is endowed with the standard  $\|.\|_\infty$-norm topology.
It is easy to see that  $\tau$ is metrizable.
It is compatible with the following metric $d_\tau$:
$$
d_\tau(T,S):=\sum_{m=1}^\infty\frac 1{2^m}\sup_{n\in\Bbb N}|\mu(T^nA_m\cap A_m)-\mu(S^nA_m\cap A_m)|,
$$
where  $(A_m)_{m=1}^\infty$ is a dense family in $\goth B$.
It is obvious that $\tau$ is stronger then the week topology.
Tikhonov shows in \cite{Ti1} that $\tau$ is not separable.
The standard action of Aut$(X,\mu)$ (endowed with the weak topology)
on (Aut$(X,\mu),\tau)$ by conjugation is continuous.
The subspace $c_0\subset l^\infty$ of sequences converging to 0 is closed
in $l^\infty$.
Moreover, it is separable and hence  Polish in the induced topology.
We note that Mix$(X,\mu)=\bigcap_{A\in\goth B}f_A^{-1}(c_0)$.
It follows that Mix$(X,\mu)$ is closed  in (Aut$(X,\mu), \tau)$ and Polish in the induced topology.
The action of Aut$(X,\mu)$ on Mix$(X,\mu)$ is {\it topologically transitive}, i.e. there is a mixing transformation $T$ whose conjugacy class is $\tau$-dense in Mix$(X,\mu)$.
Hence the set of all mixing transformations with this property is an invariant dense $G_\delta$ in Mix$(X,\mu)$.
This set contains all transformations isomorphic to Cartesian products of two non-trivial (mixing) transformations \cite{Ti1}.
Hence the conjugacy class of each Bernoullian or mixing Gaussian   transformation is dense in Mix$(X,\mu)$.
The following problem is open
\roster
\item"{\bf (Pr4)}"
Is the conjugacy class of every mixing transformation  dense in
Mix$(X,\mu)$?
\endroster

\comment

The following statement can be viewed as a partial answer to this problem.

\proclaim{Theorem 12.8 \cite{Ry}} If $T\in\text{\rom{Mix}}(X,\mu)$, $A\subset X$ with $\mu(A)>0$ and $\epsilon>0$, then there is $R\in\text{\rom{Aut}}(X,\mu)$
with $\|f_A(RTR^{-1})-f_A(\text{\rom{Id}})\|_\infty<\epsilon$.
\endproclaim
\endcomment

 It is shown in \cite{Ti1} that the subset of mixing transformations with simple spectrum is residual in
\rom{Mix}$(X,\mu)$.
The following refinement of Theorem~5.2 is perhaps the most bright application of
the topology $\tau$.

\proclaim{Theorem 6.7 \cite{Ti2}} For each $n>1$, there is a mixing $T$ such that $\Cal M(T)=\{n\}$.
\endproclaim
\demo{Idea of the proof}
We use without explanation the notation introduced in the proof of Theorem~5.2.
Let
$$
\Cal A_{G_n}^0:=\{T\in\Cal A_{G_n}\mid T_{e_{*}}\text{ is mixing}\}.
$$
This set is endowed with the weakest topology in which  the maps
$$
\align
\Cal A_{G_n}^0 &\ni T\mapsto T_{g}\in\text{\rom{Aut}}(X,\mu),\  g\in G_n, \text{ and }\\
\Cal A_{G_n}^0 &\ni T\mapsto T_{e_{*}}\in(\text{\rom{Mix}}(X,\mu),\tau)
\endalign
$$
 are all continuous.
 It is easy to see that this space is Polish.
 The subsets
 $$
 \align
 \Cal A^1 &:=\{T\in\Cal A_{G_n}^0\mid T_{e_1}\text{ has a simple spectrum}\},\\
\Cal A^2 &:=\{T\in\Cal A_{G_n}^0\mid T_{e_je_1^{-1}}\text{ has a continuous spectrum if}\ 2\le j\le n\},\\
\Cal A^3 &:=\{T\in\Cal A_{G_n}^0\mid T_{e_1e_0}\text{ has a simple spectrum}\}
\endalign
$$
are  invariant $G_\delta$ in $\Cal A_{G_n}^0$.
The most difficult part of the proof is to show that they are all dense.
In view of (Pr4), it is not enough to show that they contain a free  $G_n$-action as in the proof of Theorem~5.2.
Instead, an involved approximation argument is elaborated in \cite{Ti2}.
 Now take any $T$ from the intersection $\Cal A^1\cap\Cal A^2\cap\Cal A^3$, which is non-empty.
In view of Proposition~5.1,
the transformation $T_{e_{*}}$ is as desired.
\qed
\enddemo

\head 7. Spectral multiplicities of infinite measure preserving systems
\endhead

Suppose now that $T$ is a measure preserving transformation
of an infinite $\sigma$-finite standard measure space $(X,\goth B,\mu)$.
Since $L^2(X,\mu)$ does not contain non-trivial constants,
the Koopman operator $U_T$ is defined in the entire space $L^2(X,\mu)$.
The absence of constants helps to solve the spectral multiplicity problem {\it completely} in the infinite measure preserving case.

\proclaim{Theorem 7.1 \cite{DaRy1}}
Given $E\subset\Bbb N$,
 there is an ergodic rigid infinite measure preserving transformation $T$
such that $\Cal M(T)=E$.
\endproclaim

Before we pass to the proof of the theorem,
we draw  attention  of the reader to two specific features of infinite measure preserving dynamical systems.
First, the {\it ergodicity} of an infinite measure preserving  transformation $T$ is no longer a spectral invariant.
Indeed,  if the maximal spectral type is continuous then  there are no  invariant subsets of finite measure.
However invariant subsets of infinite measure may exist.
Secondly,  we recall that   $T$  is said to be  {\it multiply recurrent} if for each $p>0$ and a subset $A\subset X$ of positive measure there is   a positive integer $n$ such that $\mu(A\cap T^nA\cap\dots\cap T^{np}A)>0$.
While every finite measure preserving transformation  is multiply recurrent \cite{Fu}, this is not true for infinite measure preserving systems.
For  various counterexamples we refer to a survey \cite{DaSi} and references therein.
However if $T$ is rigid then $T$ is multiple recurrent.

\demo{Idea of the proof of Theorem 7.1}
Consider first a particular case when $E=\{k\}$ for some $k>1$.
Let $T$ be an infinite measure preserving rank-one transformation
such that $\exp U_T$ has a simple spectrum.
Then $U_{T^{\odot (k-1)}\times T}=(U_T)^{\odot (k-1)}\otimes U_T$.
Hence $\Cal M(T^{\odot (k-1)}\times T)=\{k\}$ by Proposition~3.4(ii), as desired (provided that $T^{\times k}$ is ergodic and rigid).

Now consider the general case.
Fix $k\in E$, $k\ne 1$.
  Then construct $K,D,H,T,\tau,Y$ in a similar way as in the proof of  Theorem~4.1.
Moreover, on infinitely many steps of the inductive cutting-and-stacking construction of $T$ put  many additional spacers to guarantee  that $T$ is {\it infinite} measure preserving.
Next consider the product transformation  $T_{Y\backslash \tau}\times T_\psi^{\odot (k-1)}$.
It is rigid by construction.
 Slightly modifying the proof of Theorem~4.1 and applying Proposition~3.4(ii), we obtain that $\Cal M(T_{Y\backslash \tau}\times T_\psi^{\odot (k-1)})=E\cup\{k\}=E$.
A standard criterium of ergodicity of cocycles
(see \cite{Sc}, \cite{GoSi})
is used
to prove ergodicity of $T_{Y\backslash \tau}\times T_\psi^{\odot (k-1)}$. \qed
\enddemo

A transformation $T$ of infinite measure space is called {\it mixing} or {\it of zero type} \cite{DaSi} if $U_T^n\to 0$ weakly.
We note that mixing  in the infinite measure does not imply ergodicity.

\proclaim{Theorem 7.2 \cite{DaRy2}}
Given  $E\subset\Bbb N$,
there is an ergodic mixing multiply recurrent infinite measure preserving transformation $T$
such that $\Cal M(T)=E$.
\endproclaim
\demo{Idea of the proof} By Theorem~7.1, there is ergodic rigid infinite measure preserving $T$ with $\Cal M(T)=E$.
Hence $T^{l_i}\to \text{Id}$ as $i\to\infty$.
We construct a rank-one mixing transformation $S$ such that
$\Cal M(T\times S)=\Cal M(T)$ and the product $T\times S$ is ergodic.
Since $S$ is mixing,  $T\times S$ is mixing too.
To construct such an $S$, the technique to force mixing (see Section~6)
is applied.
Namely, $S$ appears as  a limit of a sequence of rank-one transformations $S_n$   such that
$$
U_{S_n}^{l_{n,k}}\to\delta_n U_{S_n}\tag7-1
$$
along
a subsequence $(l_{n,k})_{k=1}^\infty$  of $(l_i)_{i=1}^\infty$,
where  $\delta_n\to 0$.
To construct $S_n$,  an appropriate class of infinite measure preserving systems is introduced.
 A rank-one transformation  is called a {\it high staircase} if it is associated with $(r_n,s_n)_{n=1}^\infty$ such that
 $s_n(i)=z_n+i$ for all $0\le i<r_n$ and some $z_n\ge 0$.
 It is shown in \cite{DaRy2} that under an infinite version of the restricted growth condition from \cite{Ad} (see also Section~6), each  high staircase  is mixing {\it independently} of the choice of $(z_n)_{n=1}^\infty$.
This independence plays the key role to select $(l_{n,k})_k$ satisfying
\thetag{7-1} inside $(l_i)_{i}$.
Now~\thetag{7-1} and Proposition~3.8(ii) imply that $U_T$ and $U_S$
are strongly disjoint and hence $\Cal M(T\times S)=\Cal M(T)=E$.
The ergodicity of $T\times S$ is {\it forced} in \cite{DaRy2} simultaneously with forcing the mixing property
of $S$.
\qed

We conclude this section with an infinite analogue of Theorem~6.2.
It is shown in \cite{DaRy2} via the   forcing of mixing technique.

\proclaim{Theorem 7.3} There is a mixing  rank-one conservative infinite measure preserving transformation $T$ such that $\exp{U_T}$ has a simple spectrum.
\endproclaim

\enddemo

\head 8. Gaussian and Poisson realizations
\endhead

We first recall definitions of the most popular dynamical systems of probabilistic origin: Gaussian and Poisson ones.
See \cite{Co-Sin}, \cite{Le-Th}, \cite{Ne}, \cite{Roy} for more information on them.
{\subhead Gaussian systems\endsubhead}
Let  $\sigma$ be a symmetric probability measure on $\Bbb T$.
Denote by $(\widehat\sigma(n))_{n\in\Bbb Z}$ the Fourier coefficients of $\sigma$.
The Abelian group  $X:=\Bbb R^\Bbb Z$ endowed with the product topology is a Polish (but non-locally compact).
Then the dual group $\widehat X$, i.e. the group of continuous  characters of $X$, is identified naturally with $\bigoplus_{n\in\Bbb Z}\Bbb R$, i.e. the group of infinite sequences $(t_n)_{n\in\Bbb Z}$ such that $t_n\ne 0$ for finitely many $n$.
The duality form is given by
$$
(y,t)\mapsto\exp\bigg(i\sum_{j\in\Bbb Z}y_jt_j\bigg),\quad y=(y_j)_j\in X,\quad t=(t_j)_j
\in\widehat X.
$$
Each probability measure $\mu$ on $X$ is determined uniquely by its {\it characteristic function}, i.e. the Fourier transform $\widehat\mu(t):=\int_{X}\exp(i\sum_{j\in\Bbb Z}y_jt_j)\,d\mu(y)$, $t\in\widehat X$.
 Let $\mu_\sigma$ denote the only probability measure on $X$
 such that
$$
\widehat{\mu_\sigma}(t)=\exp\left(-\frac 12\sum_{n,m\in\Bbb Z}\widehat\sigma(n-m)t_nt_m\right), \quad t=(t_n)_{n\in\Bbb Z}\in\widehat X.
$$
Since $\widehat{\mu_\sigma}$ is invariant under the shift on $\widehat X$, it follows that  $\mu_\sigma$ is invariant under  the dual shift, say $T$, on $X$.
The dynamical system $(X,\mu_\sigma, T)$ (and, more generally, each system isomorphic to it) is called a  {\it Gaussian} dynamical system over the base $\sigma$.
Denote by $V_\sigma$ the unitary operator acting in the Hilbert space $L^2(\Bbb T,\sigma)$ by $V_\sigma f(z):=zf(z)$.
Then $U_{T}\oplus P_0$ is unitarily equivalent to $\exp V_\sigma$.
Recall that $P_0$ is the orthogonal projection on $\Bbb C1\subset L^2(X,\mu_\sigma)$.
It follows that
\roster
\item"$\bullet$"
If $\sigma$ has an atom then $T$ is not ergodic.
\item"$\bullet$"
If $\sigma$ is non-atomic then $T$ is weakly mixing.
\item"$\bullet$"
 $T$ has a simple spectrum if and only if
$V_\sigma^{\odot j}$ has  a simple spectrum for each $j$.
\endroster
Girsanov showed in \cite{Gi} that   given an ergodic Gaussian $T$, either $T$ has a simple spectrum
or $\# \Cal M(T)=\infty$.
\footnote{Provided that $\infty\not\in\Cal M(T)$.}
He also constructed the first example of an ergodic  Gaussian $T$ with a simple spectrum \cite{Gi}.
The first mixing Gaussian system with a simple spectrum appeared in \cite{New}.

{\subhead Poissonian systems\endsubhead}
Let $X=\Bbb R$ and let $\mu$ denote the Lebesgue measure on $X$.
Denote by $\widetilde X$ the  space of Radon measures on $X$.
We equip $\widetilde X$ with the standard Borel structure $\widetilde{\goth B}$ generated by the $*$-weak topology.
For each compact subset $K\subset X$, we consider a map $N_K:\widetilde X\to\Bbb R$ given by $N_K(\omega):=\omega(K)$.
Then there is a unique probability measure $\widetilde\mu$ on $(\widetilde X,\widetilde {\goth B})$ such that
\roster
\item"$\bullet$"
$N_K$ maps $\widetilde\mu$ to the Poisson distribution with parameter $\mu(K)$, i.e.
$$
\widetilde\mu(\{\omega\mid N_K(\omega)=j\})=\frac{\mu(K)^j\exp{(-\mu(K))}}{j!}
$$
for all $K$ and integer $j\ge 0$ and
\item"$\bullet$"
if $K_1\cap K_2=\emptyset$ then the random variable $N_{K_1}$
and $N_{K_2}$ are independent.
\endroster
Let $T$ be a $\mu$-preserving homeomorphism of $X$.
Then we define a Borel map $\widetilde T:\widetilde X
\to\widetilde X$ by setting $\widetilde T\omega:=\omega\circ T$.
This map is one-to-one and  preserves $\widetilde \mu$.
It is called the {\it Poisson suspension} of $T$.
A probability preserving transformation isomorphic to the Poisson suspension of an infinite measure preserving transformation is called a {\it Poisson} transformation.
An important property of the Poisson suspensions is that the Koopman operator $U_{\widetilde T}\oplus P_0$ is unitarily equivalent to $\exp(U_T)$ \cite{Ne}.
Recall that since $\mu$ is infinite, we consider $U_T$  in the entire
space $L^2(X,\mu)$.
It follows that
\roster
\item"$\bullet$"
If $T$ has an invariant subset of  finite positive measure  then $\widetilde T$ is not ergodic.
\item"$\bullet$"
If $T$ has no invariant subsets of finite positive measure  then $\widetilde T$ is weakly mixing.
\item"$\bullet$"
 $\widetilde T$ has a simple spectrum if and only if
$U_T^{\odot j}$ has  a simple spectrum for each $j$.
\item"$\bullet$"
If $T$ is totally dissipative, i.e. there exists a subset $B\subset X$ such that $X=\bigsqcup_{n\in\Bbb Z} T^nB$, then $\widetilde T$ is Bernoullian with infinite entropy.
\item"$\bullet$" If $X=X_1\sqcup X_2$ and $X_1$, $X_2$ are $T$-invariant subsets of infinite measure then $\widetilde{T}=\widetilde {T_1}\times\widetilde{T_2}$.
\endroster

Let $\sigma$ be a symmetric probability measure of maximal spectral type of $T$.
Consider a Gaussian transformation $T_\sigma$.
Then $U_{\widetilde T}$ is unitarily equivalent to $U_{T_\sigma}$.
Thus, given a Poisson transformation, we can always find  a Gaussian transformation which is unitarily equivalent to it.
We do not know if the converse is true.

\subhead Spectral multiplicities of Gaussian and Poisson systems
\endsubhead
We consider here  the following version of (Pr1).
\roster
\item"{\bf (Pr5)}"
Given a subset $E\subset\Bbb N\cup\{\infty\}$, is there a {\it Gaussian} or {\it Poisson} weakly mixing transformation $T$ such that $\Cal M(U_T)=E$?
\endroster

In view of the  aforementioned remark, it suffices to  consider only the Poisson realizations.
It follows from the properties of Poisson suspensions that
 the direct product of a Poisson transformation with a Bernoullian transformation  is again a Poisson  transformation.
Therefore it suffices to consider only subsets $E\subset\Bbb N$.
As was noted in \cite{Roy}, there exists a Poisson transformation with a simple spectrum.
   It remains to consider the case when $E$ is  infinite.

\proclaim{Theorem 8.1 \cite{DaRy2}}
Every multiplicative  (and hence every additive) subsemigroup of $\Bbb N$ is realizable as the set of spectral multiplicities of a  mixing Poisson transformation.
\endproclaim
\demo{Proof}
Fix $p>0$.
Denote by $E$ the multiplicative subsemigroup of $\Bbb N$ generated by $p$.
Let $T$ be a rank-one mixing transformation of an infinite measure space
$(X,\goth B,\mu)$ such that $\exp U_T$ has a simple spectrum (see Theorem~7.3).
We consider the space $(X,\mu)\times (\Bbb Z/p\Bbb Z,\lambda_{\Bbb Z/p\Bbb Z})$ and set $T_p:=T\times I$.
Then
$$
\exp U_{T_p}=\exp((U_T)^{\oplus p})=\bigoplus_{n=0}^\infty
((U_T)^{\oplus p})^{\odot n}
=\bigoplus_{n=0}^\infty
(U_T^{\odot n})^{\oplus p^n}
$$
It follows that $\Cal M(\widetilde{T_p})=\{p,p^2,\dots\}$.
 Thus $E$ is Poisson realizable.

Now let, $p\ne q\in\Bbb N$.
 Denote by $E$ the multiplicative subsemigroup of $\Bbb N$ generated by $p$ and $q$.
Let $T$ and $S$ be two mixing rank-one infinite measure preserving transformations such that the unitary operator
$\exp U_T\otimes\exp U_{S}$ has a simple spectrum.
Since
$
\exp (U_{T_p}\oplus U_{S_q})=
\exp U_{T_p}\otimes\exp U_{ S_q},
$
it follows that  $\Cal M(\widetilde{T_p\sqcup S_q})=\{p,p^2,\dots\}\diamond\{q,q^2,\dots\}=E$.
 Thus $E$ is Poisson realizable.

The general case is considered in a similar way.
 \qed
\enddemo

It may seem that for each Poissonian  (or Gaussian) transformation $T$,
the set $\Cal M(T)$ is a multiplicative  subsemigroup of $\Bbb N$.
 For instance, this is claimed in the introduction to \cite{Ro3}.
The following  counterexample  is constructed in \cite{DaRy2}.

\example{Example 8.2} Let $T$ be an ergodic infinite measure preserving transformation such that $\exp(U_T)$ has a simple spectrum.
Then
$$
U_{\widetilde{T\odot T}}=\exp(U_{T\odot T})=\exp(U_T\odot U_T)=\bigoplus_{n=0}^\infty
(U_T\odot U_T)^{\odot n}.
$$
It follows from Proposition~3.4(ii) that $\Cal M((U_T\odot U_T)^{\odot n})=\{(2n)!/(2^nn!)\}$.
Since the measures of maximal spectral type of the operators $(U_T\odot U_T)^{\odot n}$, $n\in\Bbb N$, are pairwise
disjoint, we obtain that
$$
\Cal M(\widetilde{T\odot T})=\bigg\{\frac{(2n)!}{2^nn!}\,\bigg|\,n\in\Bbb N\bigg\}=\{1,3,3\cdot 5,3\cdot 5\cdot 7,\dots\}.
$$
\endexample

We say that a pair of natural numbers $(m,k)$ is {\it good} if there is a subgroup $\Gamma$ in $\goth S_m$ such that $\# \Gamma=k$.
We then set $A(m,k):=\{{(mn)!}/(k^nn!)\mid n\in\Bbb N\}$.
Generalizing Theorem~8.1 and Example 8.2 we obtain the following.

\proclaim{Proposition 8.3}
Given a sequence of good pairs $(m_1,k_1),\dots,(m_l,k_l)$, the set
$A(m_1,k_1)\diamond\cdots\diamond A(m_l,k_l)\subset\Bbb N$ is Poisson realizable.
The same is true for an infinite sequence of good pairs $(m_1,k_1),(m_2,k_2),\dots$.
\endproclaim

\head 9. Spectral multiplicities of ergodic actions of other groups
\endhead

\subhead Ergodic flows\endsubhead
Since $\Bbb R$ contains $\Bbb Z$ as a co-compact subgroup, it seems natural to export the results on spectral multiplicities for $\Bbb Z$-actions to $\Bbb R$-actions via the inducing.
The concept of {\it induced action} was introduced in \cite{Ma} (see also \cite{Zi}).
In our case, given an ergodic transformation $T$, the induced $\Bbb R$-action $W=(W_t)_{t\in\Bbb R}$ is the flow  built under the constant function $1$ over $T$.

 \proclaim{Theorem 9.1 \cite{DaLe}} \roster
 \item"\rom{(i)}"
 Two ergodic transformations $T$ and $T'$ are isomorphic if and only if the  flows $W$ and $W'$ induced by them are isomorphic.
\item"\rom{(ii)}"
  ${\Cal M}(W)={\Cal M}(T)\cup\{1\}$.
\item"\rom{(iii)}"
Thus for each $E\subset \Bbb N$ with $1\in E$, there is an ergodic flow $W$ such that ${\Cal M}(W)=E$.
\endroster
\endproclaim

However using  the method of ``induced flows''
we can not obtain weakly mixing realizations (each induced flow has a discrete component in the spectrum).
Moreover,  we can only realize subsets containing~1.
 Adapting  the cutting-and-stacking construction from Sections~2 and 4
  we can prove the following analogues of Theorems~2.1 and 4.1.

\proclaim{Theorem 9.2 \cite{DaLe}}  For each $E\subset \Bbb R$ such that $E\cap\{1,2\}\ne \emptyset$, there is a weakly mixing  flow $W$ such that ${\Cal M}(W)=E$.
Also,   ${\Cal M }(W_t)=E$ for each real $t\ne 0$.
\endproclaim

Rokhlin problem on spectral multiplicities  can be solved also for flows.

\proclaim{Theorem 9.3 \cite{DaSo}} For each $n$, there is a weakly mixing flow $W$
such that ${\Cal M} (W)=\{n\}$.
\endproclaim

The proof is based on the method of auxiliary group actions in the same way
as the proof of Theorem~5.2.
The corresponding auxiliary group is $\Bbb R\times G_n$, where $G_n$ is the auxiliary group considered in Section~5.

\subhead Other groups
\endsubhead
Let $G$ be an Abelian non-compact locally compact second countable group.
We note that the methods of realization of spectral multiplicities from the previous sections are rather specific for $\Bbb Z$-actions.
Many of them can be adjusted for $\Bbb R$-actions.
However, for general Abelian groups we encounter with new difficulties.
Our knowledge about spectral realizations in this case is rather restricted.
We hope that the $(C,F)$-construction (see the survey \cite{Da2} and references therein) being an algebraic counterpart of the geometrical cutting-and-stacking can help to construct appropriate models.

\proclaim{Theorem 9.4 \cite{DaLe}} {\sl If $G'$ is a torsion free Abelian discrete countable group then for each $E\subset\Bbb R$, if $E\cap\{1,2\}\ne\emptyset$ then there is a weakly mixing action $W$ of $G'$ with ${\Cal M}(W)=E$.}
\endproclaim

This theorem follows from the fact that if $G'\ne\Bbb Z$ then $G'$ embeds densely into $\Bbb R$.
It remains to apply Theorem 9.2.

\proclaim{Theorem 9.5 \cite{DaLe}} {\sl If $E\ni 1$ and one of the following holds
\roster
\item"\rom{(i)}"
 $G$ has a closed one-parameter subgroup,
\item"\rom{(ii)}"
 $G=D\times F$ where $D$ is a torsion free discrete Abelian group
and $F$ is an arbitrary locally compact second countable  group,
\endroster
then there is a weakly mixing action $W$ of $G$ with ${\Cal M}(W)=E$.}
\endproclaim

\proclaim{Theorem 9.6 \cite{DaSo}}  Let $G=\Bbb R^p\times G'$, where
 $G'$ contains a compact open subgroup
$G_0'$.
If one of the following is satisfied
\roster
\item"\rom{(i)}"
 $p>0$,
\item"\rom{(ii)}"
 $p=0$ but $G_0'$ is a direct summand in  $G'$
\item"\rom{(iii)}"
 $p=0$,  $G_0'$ is not a direct summand in  $G'$ but there is no
$k>0$ such that $k\cdot g'=0$ for all $g\in G'/G_0'$.
\endroster
Then for each $n>1$, there is  an action $W$ of $G$ such that
${\Cal M}(W)=\{n\}$.
\endproclaim

\proclaim{Theorem 9.7 \cite{So}}
Let $G$ be a discrete Abelian group or $\Bbb R^m$ with $m\ge 1$.
For each (finite or infinite) sequence of positive integers $p_1,p_2,\dots$, there is a weakly mixing action $W$ of $G$
such that 
${\Cal M}(W)=\{p_1\}\diamond \{p_1\}\diamond\cdots$.
\endproclaim

\head 10. Concluding remarks and open problems
\endhead

\roster
\item
Two unitary operators $U$ and $V$ in a Hilbert spaces $\Cal H_1$  and $\Cal H_2$ respectively are called {\it cyclicly equivalent} if there is a unitary operator $W:\Cal H_1\to\Cal H_2$
such that the image under $W$ of each  $U$-cyclic subspace is a $V$-cyclic subspace and vice versa.
 It is shown in \cite{Fr} that given two weakly mixing transformations $T$ and $S$, the Koopman operators $U_T$ and $U_S$ are cyclicly equivalent if and only if $\Cal M(U_T)=\Cal M(U_S)$.
\item
How to realize $\{3,4\}$ or $\{3,5\}$?
 \item (Thouvenot's question)
Which  subsets $E\ne\{1\}$  are realizable on {\it prime} transformations?
A transformation $T$ of $(X,\goth B,\mu)$ is called prime if $\goth B$ and $\{\emptyset,X\}$ are the only (up to  $\mu$-null subsets) factors, i.e. invariant sub-$\sigma$-algebras, of $T$.
All the realizations considered in this paper have non-trivial proper factors.
\item
In view of recent preprint \cite{Pr}, it is natural to ask: which subsets $E$ are realizable on transformations with absolutely continuous spectrum.
\item
Which subsets  $E\ne\{1\}$  admit
 smooth realizations?
It is shown  in \cite{BlLe} that every finite $E$ containing $1$ is realizable as  $\Cal M(T)$ for a Lebesgue measure preserving $C^\infty$-diffeomorphism $T$ of a finite dimensional torus.
\item
Are there realizations in the class of interval exchange maps?
This problem was under consideration in \cite{Os}, \cite{Ro1}, \cite{Ag2}, \cite{Ry4}.
If $T$ exchanges $n$ intervals then
$\max\Cal M(T)\le n-1$ \cite{Os}.
For each finite $E$ containing $1$, Ageev constructs in \cite{Ag2} an ergodic interval exchange transformation $T$ with $\Cal M(T)=E$.
This generalizes earlier results from \cite{Os} and \cite{Ro1}.
Ryzhikov showed in \cite{Ry4} that there is a three interval exchange transformation $T$ such that $\exp U_T$ has a simple spectrum.
\item
A transformation $T$ of $(X,\goth B,\mu)$ is said to be of rank $k$, if $k$ is the least number such that there exists an infinite  sequence of $k$ mutually disjoint  Rokhlin towers for $T$ that approximates $\goth B$.
We then write rk$(T)=k$.
It is easy to verify that rk$(T)\ge\max\Cal M(T)$ \cite{Ch}.
It is shown in \cite{KwLa} and \cite{FiKw} that given arbitrary $l\le k$, there exists
an ergodic transformation $T$ with $l=\max\Cal M(T)$ and rk$(T)=k$.
\item
Robinson in \cite{Ro4} studied the spectral multiplicities of general isometric extensions (see Section~2).
 He showed, in particular, that given an ergodic cocycle $\phi$ of $T$ with values in a compact group $K$, then  $\max\Cal M(T_\phi)$ is no less then
   the dimension of each irreducible representation of $K$.
\item
 (Ageev's problem from \cite{Ag6}) Let $G$ be an infinite countable discrete group.
 Then by \cite{Gl--We}, it has the weak Rokhlin property, i.e. there is a $G$-action $T$  whose conjugacy class is dense in $\Cal A_G$ (see Section~5 for definitions).
  Ageev shows in \cite{Ag6} that then there is a map $m_G:G\to\Bbb N\cup\{\infty\}$ such that  $m_G(g)=\max\Cal M(T_g)$ for a residual subset of actions $T\in\Cal A_G$.
If $G$ is Abelian then $m_G(g)=1$ if $g$ has an infinite order and $m_G(g)=\infty$ otherwise.
The problem is to describe $m_G$ for non-Abelian $G$.
For instance, it is shown in Section~5 that $m_{G_n}(e_1e_0)=1$ and
$m_{G_n}(e_{*})=n$.
\item
Given a weakly mixing transformation $T$, it is interesting to investigate the spectral multiplicities $\Cal M(T^n)$ for all $n>0$.
For a generic $T\in\text{Aut}(X,\mu)$, we have $\Cal M(T^n)=\{1\}$ for all $n>0$.
For each $m\ge 1$, there are  transformations $R$
with
 $\max \Cal M(R^n)=2mn$, $n>0$ (see \cite{MatNa}, \cite{Ag1} and \cite{Le1}).
A family of transformations $S$ with ``non-trivial''  sequence  $\Cal M(S^n)$
was constructed in a recent work \cite{Ry6}.
For instance, there is a weakly mixing transformation  $S$ with
$\Cal M(S^n)=\{j_n\}$ for some $j_n>0$ such that the sequence $(j_n/n)_{n=1}^\infty$ has infinitely many limit points.
\item
In a similar way, given a weakly mixing flow $W=(W_t)_{t\in\Bbb R}$,
it is  interesting to study the map $\Bbb R\ni t\mapsto\Cal M(W_t)$.
Not a lot is known about it (see Theorem~9.2).
Can it be ``highly'' non-constant?
As shown  in a recent paper \cite{LePa},  the map
$\Bbb R\ni t\mapsto\sup\Cal M(W_t)$ is of second Baire class.
This answered a question raised by Thouvenot.
\item
Is it possible to modify
the construction of mixing  realizations from  Section~6 to make it effective?
More precisely, how to choose the sequence $(N_m)_{m=1}^\infty$ from the proof of Theorem~6.2 effectively?
\item
Investigate the spectral multiplicity problem for non-Abelian non-compact locally compact  second countable groups of type {\bf I} \cite{Ki}.
Such groups include the connected nilpotent Lie groups,
the connected semisimple Lie groups and the real or complex linear algebraic groups.
A discrete countable group is of type {\bf I} if and only if it has a normal Abelian subgroup of finite index.
For the unitary representations of type {\bf I}  groups there is an analogue of the spectral theorem \cite{Ki}.
Hence given a measure preserving action $T$ of a type {\bf I} group, the
spectral invariant $\Cal M(T)$ is well defined.
\item
A transformation $T$ of $(X,\goth B,\mu)$ is called {\it approximately transitive (AT)}  if for each finite family of non-negative functions
$f_1,\dots,f_l\in L^1_+(X,\mu)$ and any $\epsilon>0$ there exist $n_1,\dots,n_s\in\Bbb Z$ and $f\in L^1_+(X,\mu)$ such that
$\inf_{\lambda_1,\dots,\lambda_l\ge 0}\|f_j-\sum_{k=1}^s\lambda_kf\circ T_{n_k}\|_1<\epsilon$ for each $j=1,\dots,l$ \cite{CoWo}.
In a similar way, one can define the AT-property for actions of arbitrary locally compact groups.
Thouvenot observed that each AT-transformation has a cyclic vector in $L^1(X,\mu)$ (for the $L^1$-analogue of the Koopman operator).
The problem is whether the AT-property implies  simplicity of spectrum of $U_T$?
Examples of non-AT systems with simple spectrum appeared in a recent paper \cite{AbLe}\footnote{Only for actions of some Abelian torsion groups. }.

\endroster

\head Appendix. Algebraic realizations
\endhead

In order to obtain ergodic systems with various spectral multiplicities
Robinson introduces in \cite{Ro1} and \cite{Ro2} a preliminary step which we call an {\it algebraic realization} of subsets of positive integers.
This preliminary step played an important role in subsequent papers \cite{Go--Li}, \cite{KwiLe}, \cite{KaLe}, \cite{Da3}, \cite{Da4}, \cite{DaRy1} on spectral multiplicities.
This appendix is devoted completely to algebraic realizations.
We need some notation.
Let $G$ be a countable Abelian group and let $v$ be a group automorphism of $G$.
Given $g\in G$, we denote by $\Cal O(g)$ the $v$-orbit of $g$,
i.e. $\Cal O(g):=\{v^i(g)\mid i\in\Bbb Z\}$.
The cardinality of
 $\Cal O(g)$ is denoted by $\#\Cal O(g)$.
We put
$$
L(G,v):=\{\#\Cal O(g)\mid g\in G, g\ne 0\}.
$$
After Robinson showed in \cite{Ro2} that every set $L(G,v)\cup\{1\}$ with $G$ finite is realizable,
a natural problem arose:
\roster
\item"---"
{\it which subsets of $\Bbb N$ can be written as
$L(G,v)$?}
\endroster
He gave a partial answer there.

\proclaim{Lemma A.1 \cite{Ro2}} Let $E$ be a finite subset of $\Bbb N$ such that
$$
\text{if $n_1,n_2\in E$ then \rm{lcm}$(n_1,n_2)\in E.$}\tag A-1
$$
Then there is a finite Abelian group $G$ and an automorphism $v$ of
$G$ such that $E=L(G,v)$.
\endproclaim

Does this lemma extend to subsets not satisfying \thetag{A-1}?
In \cite{KwiLe} the negative answer (without proof) to this question is attributed to Weiss.

\proclaim{Proposition A.2}
Let $G$ be an Abelian group and let $v$ be an automorphism of $G$.
 If $n_1,n_2\in L(G,v)$ then \rom{lcm}$(n_1,n_2)\in L(G,v)$.
\endproclaim
\demo{Proof}
Let $n_i=\#\Cal O(g_i)$ for some $g_i\in G$, $i=1,2$.
We put
 $d:=\# \Cal O(g_1+g_2)$.
Suppose first that $n_1$ and $n_2$ are coprime.
Since
 $v^{n_1n_2}(g_1+g_2)=g_1+g_2$, it follows that $d$ divides $n_1n_2$.
Therefore there are $d_1,d_2,l_1,l_2$ such that $d=d_1d_2$,  $n_1=d_1l_1$ and $n_2=d_2l_2$.
Then
$$
g_1+g_2=v^{dl_1}(g_1+g_2)=v^{n_1d_2}(g_1)+v^{n_1d_2}(g_2)=g_1+
v^{n_1d_2}(g_2).
$$
Therefore $g_2=v^{n_1d_2}(g_2)$ and hence  $n_2$ divides $n_1d_2$.
Since $n_1$ and $n_2$ are coprime, $n_2$ divides $d_2$.
Hence $n_2=d_2$.
In a similar way, $n_1=d_1$.
Therefore $d=n_1n_2$, as desired.

In the general case, let $n_1=m_1p$ and $n_2=m_2p$ for some
$p,m_1,m_2$ such that $m_1$ and $m_2$ are coprime.
Let $w:=v^p$.
It follows from the above that the $w$-orbit of $g_1+g_2$ consists of
$m_1m_2$ elements.
Hence  the $v$-orbit of $g_1+g_2$ consists of $pm_1m_2$ elements, i.e.
lcm$(n_1,n_2)\in L(G,v)$.
\qed
\enddemo

Thus, we see that the class of subsets $L(G,v)\subset\Bbb N$   is rather restrictive.
Therefore to construct spectral realizations of
 subsets not satisfying \thetag{A-1} and
 infinite subsets,
Kwiatkowsky (jr) and  Lema\'nczyk established
 in \cite{KwiLe} a refined version of Lem\-ma~A.1.
To state it we need some more notation.
Let $G$ and $v$ be as in Lemma~A.1  and let $H$ be a subgroup of $G$.
We put $L(G,v,H):=\{\#(\Cal O(h)\cap H)\mid H\ni h\ne 0\}$.

\proclaim{Lemma A.3 \cite{KwiLe}} Let $E$ be a subset of $\Bbb N$.
If $1\in E$ then there exist  an Abelian countable group $G$, a subgroup $H$ of $G$ and an automorphism $v$  of $G$ such that
$E=L(G,v,H)$.
\endproclaim

The next task was to remove  the condition $1\in E$ from the statement of Lem\-ma~A.3.
(It is interesting to note that there was no this condition in Lemma~A.1.)
This problem was risen by Katok and Lema\'nczyk in \cite{KaLe} devoted to spectral realization of subsets containing $2$ (and not containing $1$).
They considered only finite subsets $E\subset\Bbb N$ and
they needed  finite $G$  to realize
 such subsets.
In \cite{Da3}, the condition $1\in E$ was removed from the statement of Lemma~A.3 and  spectral realizations for each $E\subset \Bbb N$ with $2\in E$ were constructed.
However this did not answer the Katok-Lema\'nczyk question because the corresponding groups $G$ were always infinite, even for finite $E$.
The problem was finally settled in \cite{Da4}.

\proclaim{Lemma A.4 \cite{Da4}}
Let $E$ be a  non-empty subset of $\Bbb N$.
Then there exist a countable Abelian group $G$, an automorphism $v$ of $G$
and a subgroup $H$ in $G$ such that
$E=L(G,v,H)$.
 Moreover, if  $E$ is finite then we can choose $G$ finite.
\endproclaim
\demo{Proof}
Let $E=\{p_1,p_2,\dots,\}$ with $p_1<p_2<\cdots$.
For each $i>0$, we select a finite group $B_i$ and an automorphism $\theta_i$ of $B_i$ such that for each $0\ne b\in B_i$, the length of the $\theta_i$-orbit of $b$ is $p_i$.
We now set
$$
\align
G&:=B_1\oplus B_2^{\oplus p_1}\oplus B_3^{\oplus (p_1p_2)}\oplus\cdots,\\
H&:=\{(g_1,g_2,\dots)\in G\mid g_i=0\text{ if }i\notin\{1,2,2+p_1,2+p_1+p_1p_2,\dots\}\}.
\endalign
$$
We consider two automorphisms of $B_i^{\oplus(p_1\cdots p_{i-1})}$.
The first one, say $\sigma_i$, is generated by the cyclic permutation of the coordinates.
The second one, say $\theta_i'$, is the Cartesian product $\theta_i\times\text{id}\times\cdots\times\text{id}$.
We note that $(\theta_i'\sigma_i)^{p_1\cdots p_{i-1}}=\theta_i\times\cdots\times\theta_i$.
We now define an automorphism $v$ of $G$ by setting
$$
v:=\theta_1\times \theta_2'\sigma_2\times \theta_3'\sigma_3\times\cdots.
$$
We consider $G$ equipped with $v$ as a $\Bbb Z$-module.
Let us compute the set $L(G,v,H)$.
Take $0\ne h=(g_1,g_2,\dots,)\in H$.
Let $l$ the the maximal non-zero coordinate of $h$.
If $l>2$ then
$l=2+\sum_{i=1}^kp_1\cdots p_k$ for some $k>0$.
It is easy to see that
$$
\min\{i>1\mid\text{such that }v^ih\in H\}=p_1\cdots p_{k+1}.
$$
Moreover, $v^{p_1\cdots p_{k+1}}h=(g_1,\dots,g_{l-1}, \theta_{k+2}g_l,0,0,\dots)$.
Hence
$\#(\Cal O(h)\cap H)=p_{k+2}$.
The remaining cases where $l=1,2$ are considered in a similar way.
We thus obtain $L(G,v,H)=E$.
\qed
\enddemo

We need a  {\it compactified} version of this lemma.
Suppose that $G$ is a $D$-module for a compact Abelian group $D$.
Then we let
$$
L(G,D,H):=\{\#(\Cal O(h)\cap H)\mid H\ni h\ne 0\},
$$
where $\Cal O(h)$ denotes the orbit of $h$ under the action of $D$.
We call a compact $D$-module $K$ {\it finitary} if the subset of elements with finite $D$-orbits is dense in $K$.

\proclaim{Lemma A.5 \cite{Da4}}
Let $E$ be a non-empty subset of $\Bbb N$.
Then there exist a compact monothetic totally disconnected group $D$, a countable $D$-module $G$
and a subgroup $H$ in $G$ such that
$E=L(G,D,H)$.
 Moreover, the dual $D$-module $\widehat G$ is finitary.
If  $E$ is finite then we can choose $G$ finite.
\endproclaim
\demo{Proof}
This lemma follows  from Lemma~A.4.
The desired group $D$ is the closure of the cyclic group generated by $v$ in the group of all automorphisms of $G$ endowed with the natural Polish topology.
It is clear that $L(G,D,H)=L(G,v,H)$.
\qed
\enddemo

\enddocument